\documentclass[conf]{new-aiaa}
\def\arxivyes{1}
\usepackage[utf8]{inputenc}
\usepackage{textcomp}

\usepackage{graphicx}
\graphicspath{./figs/}

\usepackage{amsmath,amsthm}
\usepackage[version=4]{mhchem}
\usepackage{siunitx}
\usepackage{longtable,tabularx}
\usepackage{color}
\usepackage{mathtools}

\usepackage{algorithm}
\usepackage{algpseudocode}

\newcommand{\qvec}[1]{[#1]}
\newcommand{\qcross}[1]{[#1]^{\times}}
\newcommand{\qcrosssmall}[1]{\bar{#1}^{\times}}
\newcommand{\dq}[1]{\boldsymbol{\mathbf{#1}}}
\newcommand{\dqvec}[1]{[\boldsymbol{\mathbf{#1}}]}
\newcommand{\rqm}[1]{[|#1|]_R}
\newcommand{\lqm}[1]{[|#1|]_L}
\newcommand{\rdqm}[1]{[|\boldsymbol{#1}|]_R}
\newcommand{\ldqm}[1]{[|\boldsymbol{#1}|]_L}
\newcommand{\Cov}{\mathrm{Cov}}

\title{Distributed Dual Quaternion Extended Kalman Filtering for Spacecraft Pose Estimation}

\author{
Mathias Hudoba de Badyn
\footnote{Associate Professor, Department of Technology Systems, University of Oslo, Gunnar Randers vei 19, 2027, Kjeller, Norway. \texttt{mathias.hudoba@its.uio.no}, AIAA Member},
Jonas Binz
\footnote{Master's Student, Automatic Control Laboratory, Physikstrasse 3, 8092, Z\"{u}rich, Switzerland.},
Andrea Iannelli
\footnote{Assistant Professor, Institute for Systems Theory and Automatic Control, University of Stuttgart, Pfaffenwaldring 9, 70569 Stuttgart, Germany \texttt{andrea.iannelli@ist.uni-stuttgart.de}}, and 
Roy S.~Smith
\footnote{Professor, Automatic Control Laboratory, Physikstrasse 3, 8092, Z\"{u}rich, Switzerland. \texttt{rsmith@control.ee.ethz.ch}, AIAA Associate Fellow}
}
\affil{Department of Technology Systems, University of Oslo, 2027, Kjeller, Norway}
\affil{Automatic Control Laboratory, Department of Information Technology and Electrical Engineering, Swiss Federal Instutute of Technology (ETH), Z\"{u}rich, 8092, Switzerland}
\affil{Institute for Systems Theory and Automatic Control, University of Stuttgart, Pfaffenwaldring 9, 70569 Stuttgart}

\begin{document}

\maketitle

\begin{abstract}
In this paper, a distributed dual-quaternion multiplicative extended Kalman filter for the estimation of poses and velocities of individual satellites in a fleet of spacecraft {is analyzed}.
{The} proposed algorithm uses both absolute and relative pose measurements between neighbouring satellites in a network, allowing each individual satellite to estimate its own pose and that of its neighbours.
By utilizing the distributed Kalman consensus filter, a novel sensor and state-estimate fusion procedure {is proposed} that allows each satellite to improve its own state estimate by sharing data with its neighbours over a communication link.
A leader-follower approach, whereby only a subset of the satellites have access to an absolute pose measurement {is also examined}.
In this case, followers rely solely on the information provided by their neighbours, as well as relative pose measurements to those neighbours. 
The algorithm is tested extensively via numerical simulations, and it is shown that the approach provides a substantial improvement in performance over the scenario in which the satellites do not cooperate.
A case study of satellites swarming an asteroid is presented, and the performance in the leader-follower scenario is also analyzed.

\end{abstract}

\section*{Nomenclature}

{\renewcommand\arraystretch{1.0}
\noindent\begin{longtable*}{@{}l @{\quad=\quad} l@{}}
 $\mathbf{I}_{n\times k}$ & identity matrix in $\mathbb{R}^{n\times k}$\\
 $\mathbf{0}_{n\times k}$ & zero matrix in $\mathbb{R}^{n\times k}$\\
 $\mathbf{1}_n$ & vector of all ones in $\mathbb{R}^n$\\
 $I$ & inertial frame\\
 $B,B_i$ & body frame, body frame of satellite $i$\\
 $A_i$ & body frame of the IMU of satellite $i$\\
 $\mathbb{H},\mathbb{H}^u$ & set of quaternions/unit quaternions $\{q\in\mathbb{H}~|~qq^*=1\}$\\
 $A(q)$ & rotation matrix corresponding to quaternion $q$\\
 $q_s,\bar{q}$ & scalar and vector parts of quaternion $q$\\
 $q_{B/I}$ & unit quaternion defining the attitude transformation $I$ to $B$\\
 $r_{B/I}^I$ & vector quaternion defining the position of $B$ in respect to $I$ expressed in $I$\\
 $\overline{r_{B/I}^I}$ & three-dimensional vector defining the position of $B$ in respect to $I$ expressed in $I$\\
 $\epsilon$ & a dual number\\
 $\mathbb{H}_d,\mathbb{H}^u_d$ & set of dual quaternions, and unit dual quaternions $\{\dq{q}\in\mathbb{H}_d~|~\dq{q}\dq{q}^*=1\}$\\
 $\dq{q},q_r,q_d$ & dual quaternion, real/dual parts of dual quaternion\\
 $\dq{q}_{B/I},\overline{\dq{q}_{B/I}}$ & unit, and reduced unit dual quaternion defining the pose transformation from $I$ to  $B$\\
 $\dq{\omega}_{B/I}^B,\overline{\dq{\omega}_{B/I}^B}$ & vector and reduced dual velocity of  $B$ compared to $I$ expressed in $B$ \\
 $\omega_{B/I}^B$ & angular velocity of $B$ compared to $I$ expressed in $B$ [rad/s]\\
$v_{B/I}^B$ & linear velocity of $B$ compared to $I$ expressed in $B$ [m/s]\\
 $[|q|]_R,[|q|]_L$ & function from $\mathbb{R}^4 \rightarrow \mathbb{R}^{4\times 4}$ denoting the right/left quaternion multiplication\\
 $[|\dq{q}|]_R,[|\dq{q}|]_L$ & function from $\mathbb{R}^8 \rightarrow \mathbb{R}^{8\times 8}$ denoting the right/left dual quaternion multiplication\\
 $x_i$ & the state of satellite $i$\\
 $\Delta x_i$ &  Kalman state update of satellite $i$ \\
 $P_i$ & the estimation covariance matrix of satellite $i$\\
 $F_i$ & state transition matrix of satellite $i$\\
 $G_i$ & process noise matrix of satellite $i$\\
 $z_i$ & measurement vector of satellite $i$\\
 $z_{i,k}$ & measurement vector of satellite $i$ expressed in the states of satellite $k$\\
 $H_i$ & observation matrix of satellite $i$\\
 $H_{i,k}$ & observation matrix of satellite $i$ expressed in the states of satellite $k$\\
 $\dq{\hat q}_{B/I},\dq{q}_{B/I, m}$ & estimate/measurement of the pose transformation from I to B\\
 $\dq{\hat \omega}_{B/I}^B,\dq{\omega}_{B/I, m}^B$ & estimate/measurement of the dual velocity of  B compared to I expressed in B \\
 $\dq{b}_{\omega},\dq{\hat b}_{\omega}$ & dual bias/estimate of dual bias\\
 $\mathbf{J}_{8\times 8}$ & extended inertia matrix\\
 $\mathcal{G}$ & communication graph of the satellite network\\
 $V$ & indexed set of nodes of $\mathcal{G}$ corresponding to individual satellites\\
 $L,F$ & sets of leader/follower nodes\\
 $E$ & set of edges of $\mathcal{G}$ corresponding to commmunication links between satellites \\
 $N_i$ & set of all neighbours of satellite $i$, excluding $i$\\
 $J_i$ & set of all neighbours of satellite $i$, including $i$\\
 $\Lambda_{i, k}$ & set of all neighbours which satellite $i$ and satellite $k$ have in common\\
 $\overrightarrow{\hat q_i}$ & estimate of all satellite poses $\in J_i$ done by satellite $i$\\
 $\overrightarrow{\hat b_i}$ & estimate of all dual bias' $\in J_i$ done by satellite $i$\\
 $\overrightarrow{\hat \omega_i}$ & estimate of all dual velocities $\in J_i$ done by satellite $i$\\
 $\overrightarrow{\hat \omega}_{i, m}$ & dual velocities measurements $\in J_i$ done by satellite $i$\\
 $y_i$ & aggregated measurement vector of satellite $i$\\
 $S_i$ & aggregated covariance matrix of satellite $i$\\
 $\mu_{q, i}$ & weight of the soft consensus of satellite $i$ on the quaternion\\
 $\mu_{r, i}$ & weight of the soft consensus of satellite $i$ on the position\\
 $\mu_{b, i}$ & weight of the soft consensus of satellite $i$ on the bias\\
\end{longtable*}}

\section{Introduction}

Terrestrial constraints such as light pollution, the day-night cycle, and an increasingly more crowded low-Earth orbit all point towards the future of astronomy being on spaceborne satellites in deep space~\cite{orbital_debris}.
Distributing the task of observation over multiple satellites has the potential to increase observational capacity by extending the baseline for imaging~\cite{bandyopadhyay2017optimal}, or providing a disturbance-free environment for novel telescope types like laser interferometers~\cite{danzmann1996lisa,amaro2017laser}.
Such satellite systems require the ability for every satellite in the fleet to point at the same object cooperatively, or maintain some relative formation, including both position and attitude (in other words, \emph{pose}) constraints.
This requires the development of novel distributed pose control and estimation algorithms.

Attitude control for spacecraft using quaternions is a mature field~\cite{Attitude_Error_Representations, FARRELL1970419, Survey_of_Nonlinear_Attitude_Estimation_Methods, 7171823}.
Conversely, the paradigm of dual quaternions, which allows a compact representation of attitude and position (or the pose) of a spacecraft, is a new area of research.
{Novel control scenarios using dual quaternions include distributed control methods~\cite{8795869}, powered-descent guidance~\cite{lee2017constrained}, data-driven techniques~\cite{koopmandq}, as well as static distributed estimation~\cite{7039403}.
Dual quaternions also provide a natural formalism for multi-spacecraft proximity operations such as rendezvous~\cite{DONG201687,filipe2015adaptive,yang2019adaptive} and multi-spacecraft formation flight~\cite{coupled_6-DOF_control}.}

Much effort has been dedicated to spacecraft estimation~\cite{Survey_of_Nonlinear_Attitude_Estimation_Methods,optimal_estimation_of_dynamic_syxstems}. 
Kalman filtering based on the unit quaternion is well-understood, and has been used in many different applications, including NASA spacecraft~\cite{Survey_of_Nonlinear_Attitude_Estimation_Methods, FARRELL1970419, 1643403}. 
The success of the quaternion multiplicative extended Kalman filter (Q-MEKF) can be attributed to different key components of the filter. 
In contrast to Euler angles, unit quaternions provide a non-singular representations of attitude, while using the minimum number of parameters\footnote{It should be noted that the Q-MEKF, which only uses the vector part of the quaternion, can still be singular if the error between the attitude and its estimate is greater than 180$^\circ$}. 
Additionally, updating attitudes in the unit quaternion framework has a lower computational complexity in comparison to Euler angles.

Dual quaternion multiplicative extended Kalman filtering (DQ-MEKF), on the other hand, is a more recent development and has not been used to the same degree as the Q-MEKF \cite{7171823,filipe2015extendedjgcd}. 
However, the compact representation of the pose via the dual quaternion make their use in spacecraft pose estimation appealing for several reasons. 
Firstly, the kinematic and dynamic equations in the dual quaternion setting have the same form as that in the attitude-only quaternion setting, meaning that much of intuition on quaternion-based estimation and control can be re-interpreted in the dual quaternion algebra.
Secondly, visual navigation systems, such as Position Sensing Diodes~\cite{visnav2, alonso2000vision}, or camera-based relative measurements~\cite{sung2022optical,filipe2015extendedjgcd, mourikis2009vision}, inherently couple attitude and translation.
Therefore, a unified framework for pose estimation is useful when considering these types of sensors.

Distributed Kalman filtering has also been an area of rapid recent development~\cite{dkf_review}. 
A distributed estimation procedure based on dual quaternions was proposed in~\cite{7039403}, wherein each sensor node estimates the pose through the Kalman filter and then updates its estimate by solving a minimization problem which is represented by a function of the estimates and covariances of its neighbours. 
A DQ-MEKF was used in a two-satellite setting for pose estimation and control in~\cite{zivan2018dual}.
A more general approach is that of the Kalman-consensus filter, where a consensus term is added after the update step of the Kalman filter~\cite{1583486, 1583238, 4434303}. 
This approach also assumes a multi-sensor model in which all sensor nodes estimate the same state --- other works with similar assumptions include~\cite{ryu2019distributed, 1583486, LI2019104500}. 
A multi-robot distributed Kalman filter where each robot estimates only its own state (in 2D), but uses relative measurements between neighbouring robots, is proposed by~\cite{1067998}. 
A clear gap in the literature is a distributed version multiplicative extended Kalman filter in the dual quaternion framework. 
We address this here, for general network topologies, and for several sensor models.

The contributions of this paper are as follows.
First, we present a novel distributed dual quaternion extended Kalman filter operating on a fleet of satellites defined on a network described by a graph.
Each satellite has access to a pose measurement, but can also measure the pose of its neighbours relative to its own pose.
Facilitated by communication between neighbouring satellites, measurement, covariance and estimation data can be exchanged, by which the pose estimates can be further improved by sensor and covariance fusion.
In particular, we propose two sensor and state estimate fusion steps, which we dub `soft' and `hard' consensus, which allow neighbouring satellites to share information about their state estimates and relative measurements and improve the overall performance of each satellite in the fleet.
Briefly, soft consensus is a distributed averaging of the estimates with neighbouring satellites, and hard consensus includes a measurement and covariance fusion step.
The proposed distributed algorithm is then tested in numerical simulations and compared to the setting in which each satellite computes its own state estimate without any cooperation. 
The simulations show that the proposed algorithm results in smaller Root Mean Square (RMS) errors in pose and velocity compared to a satellite fleet without any cooperation.

Second, we examine the filter algorithm in the case where there is a set of `leaders', which have access to a `global' pose measurement with respect to an inertial frame (e.g., a star-tracker for attitude and a laser range-finder/GPS to a nearby body for position), and a set of followers which have only access to relative measurements with respect to their neighbours.
We derive the necessary algebraic steps for filtering, sensor and measurement fusion, and communication.
Numerical experiments show that we are able to remove the absolute pose measurement from a substantial part of the satellite fleet without degradation in performance of the estimator.
Third, we provide an open-source repository~\cite{binzcode} containing the code to reproduce the numerical examples in the paper, which can also be used for more general distributed filtering problems.

The paper is organized as follows.
We introduce preliminaries in \S\ref{sec:preliminaries}, and we derive our distributed filter, including the hard and soft consensus terms, and the leader-follower setting in \S\ref{sec:distr-dual-quat}.
Numerical experiments are presented in \S\ref{sec:numerical-examples}, and the paper is concluded in \S\ref{sec:conclusion}.
\if\arxivyes1
Lengthy calculations and a primer on the DQ-MEKF are delegated to appendices.
\else
Lengthy calculations are delegated to appendices.
\fi

\section{Preliminaries}
\label{sec:preliminaries}

Our notation follows that of~\cite{filipe2015extendedjgcd}.
Additionally, for matrices $A_1,\dots,A_n$ we define $\mathrm{blkdiag}\{A_1,\dots,A_n\}$  to be the matrix formed by placing $A_1,\dots,A_n$ on the block-diagonal.
\if\arxivyes1
A summary of the DQ-MEKF from~\cite{filipe2015extendedjgcd} is provided in Appendix~\ref{sec:dq-mekf}.
\else
A summary of the DQ-MEKF from~\cite{filipe2015extendedjgcd} is provided in Appendix B of supplementary material~\cite{arxivversion}.
\fi

\subsection{Information Form Extended Kalman Filter}
\label{sec:inform-form-extend}
Consider the nonlinear dynamical system,
\begin{align}
\dot x (t) = f(x(t), t, u(t)) + g(x(t), t)w(t)\label{eq:def_kalman_state_eq}
\end{align}
with state variable $x(t)$ $\in$ $\mathbb{R}^n$, input $u(t) \in \mathbb{R}^q$, and Gaussian process noise $w(t)$ $\in$ $\mathbb{R}^p$ with $E[w(t)] = \mathbf{0}_p$ and $\Cov[w(t)] = Q(t)$.
The function $f: \mathbb{R}^n\times\mathbb{R}\times\mathbb{R}^q\mapsto \mathbb{R}^n$ maps the state and input to its derivative, and $g: \mathbb{R}^{n} \times\mathbb{R} \mapsto \mathbb{R}^{n\times p}$  maps the noise onto the dynamics.

Consider an estimate of $x(t)$ defined as $\hat{x}(t):= E[x(t)]$, and the covariance of that estimate $P(t):= \Cov[x(t)]$.
The estimate and its covariance satisfy respectively,
\begin{align}
  \dot{\hat{x}}(t) &= E[f(x(t), t, u(t))] \label{eq:estimate_de}\\
  \dot P(t) &= F(t) P(t) + P(t) F^T(t) + G(t)Q(t)G^T(t),\label{eq:Ricatti}
\end{align}
where 
\begin{subequations}\label{eq:F_and_G_def}
  \begin{align}
F(t) &= \left.\frac{\partial f(x, t, u(t))}{\partial x}\right|_{\hat x(t)}\\
G(t) &= g(\hat x(t), t).    
  \end{align}
\end{subequations}
Eq.~\eqref{eq:estimate_de} is usually approximated as,
\begin{equation}\label{eq:kalman_time_update}
    \dot{\hat{x}}(t) \approx f(\hat x(t), t, u(t)).
\end{equation}

Consider the measurement $z(t_k) \in\mathbb{R}^m$ given by the model,
\begin{equation}\label{eq:def_measurement_relationship_kf}
z(t_k) = h(x(t)) + v(t),
\end{equation}
where $v(t)$ is measurement noise represented as a Gaussian process with $E[v(t)] = \mathbf{0}_m$ and $\Cov[v(t)] = R(t)$. 

We can write the discrete measurement update as,
\begin{align}
\hat x^+(t_k) = \hat x^-(t_k) + \Delta\hat x(t_k)
\label{eq:Kalman_output}
\end{align}
where $\hat x^-$ is the predicted value given by the time update from the previous timestep, $\hat x^+$ is the predicted value after the measurement update, and $\Delta\hat x(t_k)$ is the  Kalman state update given by,
\begin{equation}\label{eq:Kalman_state_update}
    \Delta\hat x(t_k) = K(t_k)[z(t_k) - \hat z(t_k)],
\end{equation}
where {the predicted measurement is,}
\begin{align}
  \hat z(t_k) = h(\hat{x}(t)),\label{eq:meas_jacobian}
\end{align}
and the Kalman gain $K\in\mathbb{R}^{m\times n}$ is given by
\begin{align}\label{eq:Kalman_Gain}
    K(t_k) &= P^-(t_k)H^T(t_k)\left[H(t_k)P^-(t_k)H^T(t_k) + R\right]^{-1},
\end{align}
where $H(t_k)$ is the Jacobian of the measurement function,
\begin{equation}\label{eq:Kalman_H_Matrix}
    H(t_k) = \left.\frac{\partial h(x)}{\partial x}\right|_{\hat x^-(t_k)}.
\end{equation}

After the measurement, the covariance is updated by
\begin{equation}\label{eq:Kalman_H_Matrix}
    P^+(t_k) = (\mathbf{I} - K(t_k)H(t_k))P^-(t_k).
\end{equation}
Finally, $\hat{x}^+(t_k)$ is propagated forward to $t_{k+1}$ via the dynamics~\eqref{eq:kalman_time_update}, yielding $\hat{x}^-(t_{k+1})$.

Following \cite{5399678}, we introduce the information form Kalman filter, modified to the current setting of an EKF.
The information form Kalman filter has some benefits in comparison to the usual Kalman filter in the distributed case; in particular it allows batch measurement processing with lower complexity by avoiding several superfluous matrix inversions. 
From Eq.~\eqref{eq:def_kalman_state_eq} to Eq.~\eqref{eq:Kalman_output} the steps of the filter are the same. 
We introduce:
\begin{align}\label{eq:y_S_inf_form_def}
    u(t_k) &= H^T(t_k)  R^{-1} z(t_k) \\ U(t_k) &= H^T(t_k)  R^{-1} H(t_k)
\end{align}
\begin{equation}\label{eq:M_inf_form_def}
    M(t_k) = \left(P^{-1}(t_k) + U(t_k)\right)^{-1}.
\end{equation}
Finally, the  Kalman state update is computed by
\begin{equation}
    \Delta\hat x(t_k) = M(t_k)\left(u(t_k) - U(t_k)\hat x(t_k)\right).
\end{equation}
The  matrix $M$ can be seen as the \emph{a posteriori} covariance.


\subsection{Dual Quaternions}
For reasons of brevity, we assume that the reader is familiar with the basic properties and operations of quaternions, as in~\cite{optimal_estimation_of_dynamic_syxstems}.
These operations have algebraic representations that are useful for computation; in this setting, the quaternion is represented as a vector $q\in\mathbb{R}^4$. 
In order to multiply two quaternions $a$ and $b$, we define the following
\begin{align}
  \qvec{a} =   
  \begin{bmatrix}
   a_0&
   a_1&
   a_2&
   a_3
  \end{bmatrix}^T := 
  \begin{bmatrix}
    a_0 \\ \bar{a}
  \end{bmatrix}
,~
\qcross{a} := 
  \begin{bmatrix}
   0 & \mathbf{0}_{1\times 3} \\
   \mathbf{0}_{3 \times 1} & \qcrosssmall{a} 
  \end{bmatrix},~
  \qcrosssmall{a}  := 
    \begin{bmatrix}
   0 & -a_3 & a_2 \\
  a_3 & 0 & -a_1\\
  -a_2 & a_1 & 0 
  \end{bmatrix},
\end{align}
where we have adopted the scalar-first convention.
We can respectively write the \emph{left and right quaternion multiplication operators} as,
\begin{align}
\lqm{a} :=     \begin{bmatrix}
   a_0 & -\bar a^T \\
  \bar a & [\tilde{{a}}]_L\\
  \end{bmatrix},~ \rqm{a} :=     \begin{bmatrix}
   a_0 & -\bar a^T \\
  \bar a & [\tilde{{a}}]_R\\
  \end{bmatrix} \label{eq:quat_vec_def}
\end{align}
where
$[\tilde{{a}}]_L = a_0\mathbf{I}_{3\times 3} + \qcrosssmall{a}$ and $[\tilde{{a}}]_R = a_0\mathbf{I}_{3\times 3} - \qcrosssmall{a}$.
Then, quaternion multiplication and the cross product are given by
$\qvec{ab} = \lqm{a}\qvec{b}  = \rqm{b}\qvec{a}\in \mathbb{R}^{4}$
and
$\qvec{a\times b} = \qcross{a} \qvec{b} \in \mathbb{R}^4$, respectively.

The unit quaternions are used for attitude estimation as they describes a rotation from one coordinate frame  (for example, the inertial frame $I$) to another coordinate frame (for example, the body frame $B$). 
Here, $q^*$ denotes the quaternion conjugate.
Formally, a rotation amplitude $\theta$ and a rotation axis $\overrightarrow{n} \in \mathbb{R}^3$, which describes the rotation from frame $I$ to $B$, can be described as the quaternion $q_{B/I}$
\begin{equation}\label{eq:rot_quat_def}
    q_{B/I} = [\cos({\theta}/{2}), \overrightarrow{n}^T\sin({\theta}/{2})]^T.
\end{equation}
The inverse rotation is then simply $q_{B/I}^*$. 
The scalar part of the quaternion can be computed as,
\begin{equation}\label{eq:unit_quat_scal_calc}
    q_0 = \sqrt{1- \|\bar{q}\|^2}.
\end{equation}

Dual quaternions are an extension of quaternions that encode both the position and attitude of a rigid body. 
The \emph{dual quaternion} is defined by,
  \begin{align}
    &\dq{q} = q_r + \epsilon q_d \qquad \qquad q_r, q_d \in \mathbb{H},\label{eq:dquat_def}
  \end{align}
where the \emph{dual number} $\epsilon$ satisfies $\epsilon^2 = 0$ and $\epsilon \neq 0$.
The basic operations on dual quaternions are analogous to that of the quaternions, and are as follows:\\
\emph{Addition:}
\begin{equation}\label{eq:dquat_add_def}
    \dq{a} + \dq b = ({a_r} + {b_r}) + \epsilon({a_d} + {b_d}) \in \mathbb{H}_d
\end{equation}
\emph{Multiplication by scalar:}
\begin{equation}\label{eq:dquat_scalar_def}
    \lambda\dq{a} = (\lambda{a_r})+\epsilon(\lambda{a_d}) \in \mathbb{H}_d
\end{equation}
\emph{Multiplication:}
\begin{equation}\label{eq:dquat_mult_def}
    \dq{ab} = (a_rb_r) + \epsilon(a_rb_d + a_db_r) \in \mathbb{H}_d
\end{equation}
\emph{Conjugation:}
\begin{equation}\label{eq:dquat_conj_def}
    \dq{a}^* = {a^*_r} + \epsilon{a^*_d} \in \mathbb{H}_d
\end{equation}
\emph{Cross Product:}
\begin{equation}\label{eq:dquat_cross_def}
    \dq{a}\times\dq{b} = {a_r}\times{b_r}+\epsilon({a_d}\times{b_r}+{a_r}\times{b_d}) \in \mathbb{H}_d.
\end{equation}
As in the case of the quaternion, it is convenient to represent the dual quaternion operations in algebraic terms.
Define the following vectors and operators:
\begin{equation}\label{eq:quat_vec_def}
    \begin{gathered}
\dqvec{q} = \begin{bmatrix}
   \qvec{q_r} \\
   \qvec{q_d} 
  \end{bmatrix} \in \mathbb{R}^8 \qquad
    \rdqm{a} = \begin{bmatrix}
   \rqm{a_r} & \mathbf{0}_{4\times 4} \\
   \rqm{a_d} & \rqm{a_r} 
  \end{bmatrix} \in \mathbb{R}^{8\times 8}\\
  \ldqm{a} = \begin{bmatrix}
   \lqm{a_r} & \mathbf{0}_{4\times 4} \\
   \lqm{a_d} & \lqm{a_r} 
  \end{bmatrix} \in \mathbb{R}^{8\times 8}\qquad
\bar{\dq{a}} = \begin{bmatrix}
   \bar{a}_r \\
   \bar{a}_d 
  \end{bmatrix} \in \mathbb{R}^{6} \\
  [\dq{\tilde a}]=  \begin{bmatrix}
   [\tilde{a}_r] & \mathbf{0}_{3\times 3}\\
   [\tilde{a}_d] & [\tilde{a}_r] 
  \end{bmatrix}\in \mathbb{R}^{6\times 6}\qquad
  \bar{\dq{a}}^{\times} = \begin{bmatrix}
   \bar{a}_r^{\times} & \mathbf{0}_{3\times 3}\\
   \bar{a}_d^{\times} & \bar{a}_r^{\times} 
  \end{bmatrix}\in \mathbb{R}^{6\times 6}.
    \end{gathered}
\end{equation}
Then, we can represent dual quaternion multiplication as
$    \dqvec{ab} = \ldqm{a}\dqvec{b} = \rdqm{b}\dqvec{a} \in \mathbb{R}^8.$

The set of unit dual quaternions is of particular interest in the setting of \emph{pose} estimation and control, where the pose is a representation of both position and attitude. 

Consider an object whose $B$ frame is at $\overline{r_{B/I}^I}$ with respect to the $I$ frame, expressed in the $I$ frame. 
Its rotation from the origin to its body frame is described by an angle $\theta$ and its rotation axis, $\overrightarrow n$. 
The (real) rotation quaternion $q_r$ is computed using Eq.~\eqref{eq:rot_quat_def}. 
Then, the unit dual quaternion describing the pose of the object is
\begin{align}
    q_d = \frac{1}{2}r_{B/I}^Iq_r, \quad
    r_{B/I}^I = (0, \overline{r_{B/I}^I}),\quad
    \dq{q} = q_r + \epsilon q_d  \in \mathbb{H}_d^u    
\label{eq:pos_quat_def}
\end{align}

\subsection{Rigid Body Model in Dual Quaternions}
We now introduce the kinematics and dynamics for a single rigid body using dual quaternions~\cite{filipe2015extendedjgcd,coupled_6-DOF_control}. 
\subsubsection{Dual Kinematics}
The unit dual quaternion which defines the transformation from inertial frame $I$ to the body frame $B$ of the rigid body is
\begin{equation}
    \dq{q}_{B/I} = q_r + \epsilon q_d = q_r + \epsilon \frac{1}{2}r_{B/I}^Iq_r = q_r + \epsilon \frac{1}{2}q_rr_{B/I}^B.
\end{equation} 
The dual velocity $\dq{\omega}_{B/I}^B \in \mathbb{H}^v_d$, which is a vector quaternion, is defined as
\begin{equation}
\label{eq:dual_velocity_def}
\dq{\omega}_{B/I}^B = 
\begin{bmatrix} 
0 \\ {\omega}_{B/I}^B 
\end{bmatrix} + \epsilon 
\begin{bmatrix} 0 \\ {v}_{B/I}^B 
\end{bmatrix} 
\end{equation}
where $\dq{\omega}_{X/Y}^Z$ describes the \emph{dual velocity} in the $X$ frame compared to the $Y$ frame, expressed in the $Z$ frame. 
The real part describes the angular velocity, whereas the dual part describes the linear velocity.
Using the dual velocity, the dual kinematics are derived using Eq.~\eqref{eq:dquat_def} as,
{
\begin{align}
\dot{\dq{q}}_{B/I} &= \dot q_r + \epsilon \dot q_d \label{eq:derivation_kinematics}\\
                  &= \frac{1}{2}\dq{q}_{B/I}\dq{\omega}_{B/I}^B,
\end{align}
}
where the form of $\dot{q}_d$ is derived as,
\begin{align}                
\dot q_d &= \frac{1}{2} v_{B/I}^I q_r + \frac{1}{2}r_{B/I}^I \dot q_r\label{eq:derivation_qd_dot}\\
         &= \frac{1}{2} q_r v_{B/I}^B  + \frac{1}{2} q_d \omega_{B/I}^B.
\end{align}

\subsubsection{Dual Dynamics}
\begin{figure}
\begin{center}
\includegraphics[width=0.3\textwidth]{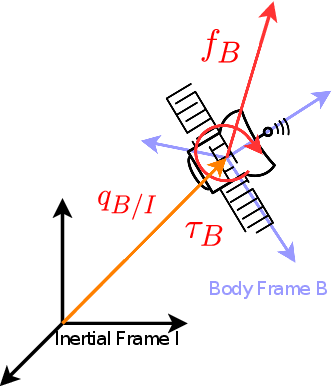}
\caption{The dual quaternion $\dq{q}_{B/I}$ defines the pose of the body frame $B$ with respect to the intertial frame $I$. Forces $f_B$ and torques $\tau_B$ are applied in the body frame.}\label{fig:rigid_body_model}
\end{center}
\end{figure}

All forces and torques from actuators are exerted on the satellite in the body frame; thus the dual force $\dq{f}_B$ is introduced as,
\begin{equation}
\label{eq:definiton_dual_force}
\dq{f}_B = 
\begin{bmatrix} 
0 \\ {\tau}_B 
\end{bmatrix} + \epsilon 
\begin{bmatrix} 
0 \\ {f}_B 
\end{bmatrix}.
\end{equation}
An illustration is shown in Fig.~\ref{fig:rigid_body_model}.
The extended inertia matrix is defined as
\begin{align}
  \label{eq:definition_Jacobian}
  \mathbf{J}_{8\times 8} = \mathrm{blkdiag}
  \left\{1,m\mathbf{I}_{3\times 3},1,\mathbf{J} \right\},
\end{align}
  where $m$ is the mass of the object and $\mathbf{J}$  is its inertia matrix. 
Following \cite{coupled_6-DOF_control}, the dual dynamics are given by,
\begin{equation}\label{eq:background_dual_dynamics}
\dq{\dot \omega}_{B/I}^B = \mathbf{J}^{-1}(\dq{f}_B - \dq{\omega}_{B/I}^B \times \mathbf{J} \dq{\omega}_{B/I}^B).
\end{equation}

\subsection{Graph Theory and Distributed Kalman Filtering}
\label{sec:distributed_KF}

When performing Kalman filtering across a network of distributed computational nodes with access to only local information, a centralized algorithm for state estimation is infeasible.
Furthermore, each node may only have access to a measurement of a subset of the states of the system.
These were some of the motivating factors for introducing distributed Kalman filtering for sensor networks in~\cite{1583238}. 
This is the foundation for the distributed approach of our present work and is therefore explained in detail here.

A graph describes the connectivity of a set of objects. 
In particular, a graph $G = (V,E)$ is a set of \emph{nodes} $V = \{1,\dots,l\}$ and a set of \emph{edges} $E\subseteq V\times V$, where  for $i,j\in V$, $ij\in E$ if and only if there is a connection between $i$ and $j$.
In the context of this paper, a `connection' between satellites $i$ and $j$ denotes that satellite $i$ and $j$ can communicate, as well as measure relative poses of each other.
For simplicity, we assume that $G$ is \emph{undirected}, in that if $ij\in E$, then $ji \in E$: if satellite $i$ can communicate with $j$, then $j$ can also communicate with $i$\footnote{If we were to consider directed graphs, then this would imply one-way communication and measurements between neighbouring satellites}.
The \emph{neighbourhood} $N_i$ of $i\in V$ is defined as the set of $j\in V$ such that $ij\in E$.
We later use a set $V_i:=N_i\cup \{i\}$ to describe the neighbourhood of $i$ if $i$ is defined to be a neighbour of itself.

Let us assume a graph $G$ with connectivity defined by an adjacency matrix $\mathcal{A}$, where each node $i$ has a set of neighbours $N_i$. 
Furthermore, the nodes want to estimate the state $x$ of the dynamical system given by Eq.~\eqref{eq:def_kalman_state_eq}\footnote{In our work, each node will only estimate a subset of the states of the system. This will be elaborated on in \S\ref{sec:derivation}.}.

Every node has an estimate $\hat x_i$ of the state $x$\footnote{Later in the paper, we adapt the distributed Kalman filter so that each satellite $i$ estimates its own state, and that of its neighbours}. 
As such, the time update can be executed locally, using the EKF procedure described in \S\ref{sec:inform-form-extend}. 
Every node propagates its estimate and covariance locally. 
However, in such a setting a node may not be able to observe all states, or each node might have a different sensor. 
Even if all nodes had a full state measurement, it may be advantageous to fuse this information to get better estimates than one node might have by itself.

Let each node $i\in V$ have the sensor model
\begin{equation}\label{eq:distributed_sensor_model}
z_i(t_k) = h(x(t_k)) + v_i(t_k).
\end{equation}
For `local' (read: non-distributed) Kalman filtering, only this measurement  would be used to compute the  Kalman update.
In the distributed setting~\cite{1583238}, measurement, estimate and covariance data from the neighbours of $i$ can be fused to improve the estimate $\hat{x}_i$ of $i$.

At each timestep $t_k$, each node $i$ first does the time update of its estimate $\hat x_i$ and covariance $M_i$ by using Eqns.~\eqref{eq:Ricatti}, \eqref{eq:F_and_G_def} and \eqref{eq:kalman_time_update}. 
Then, each node $i\in V$ takes a measurement $z_i$ with the measurement model given in Eq.~\eqref{eq:distributed_sensor_model}. 
Each node can now compute $u_i$ and $U_i$ with Eq.~\eqref{eq:y_S_inf_form_def}, which is then broadcast along with the state estimate to its neighbours.
Explicitly, the data the neighbours receive is the packet,
\begin{equation}\label{eq:distributed_message}
\text{msg}_i = (u_i, U_i, \hat{x}^-_i).
\end{equation}
Now the data and covariance matrices are aggregated in a data fusion step as,
\begin{equation}
  y_i = \sum\limits_{j \in N_i} u_j\qquad
  S_i = \sum\limits_{j \in N_i} U_j\qquad
  M_i = \left(P_i^{-1} + S_i\right)^{-1}\label{eq:M_matrix}
\end{equation}
and the Kalman state update is given by,
\begin{align}
    \hat{x}^+_i = \hat{x}^-_i + M_i[y_i - S_i\hat{x}^-_i] + \mu_i\sum\limits_{j \in N_i}\left(\hat{x}^-_j - \hat{x}^-_i\right)\label{eq:olfati_state_update}
\end{align}
where $\mu_i>0$ is a tuning parameter.
{The assumptions surrounding the optimality of the state update in Eq.~\eqref{eq:olfati_state_update} and are discussed in~\cite{5399678}.}

\section{Distributed Dual Quaternion Extended Kalman Filter}
\label{sec:distr-dual-quat}

In this section, we present the distributed dual quaternion extended Kalman filter for spacecraft pose estimation. 
The derivation of the distributed filter is conducted by assuming that each satellite has access to absolute pose measurements, after which we elaborate on the \emph{leader-follower} scenario where only some satellites have absolute pose measurements.

\subsection{Problem Statement}\label{sec:problem_statement}
Let us assume a fleet of $l$ satellites, whose communication topology is described by a graph $\mathcal{G} = ({V},{E})$.
Each node $i\in {V}$ represents a satellite, which has a set of neighbours $N_i$ with $|N_i| = k_i$. 
We emphasize that ${V}$ is an ordered set, i.e., each satellite has a numerical label $i\in {V}$, and we assume that the fleet knows its own label and those of its neighbours.

We assume that at each timestep, every satellite $i$ takes $k_i$ relative pose measurements in its own body frame, one for each neighbouring satellite.
For now, we will also assume that each satellite takes one absolute pose measurement with respect to the inertial frame, and we consider cases where the satellites may or may not have access to its velocities and accelerations in the body frame. 
We will relax this assumption later in the paper.

\begin{figure*}
\centering
\includegraphics[width=.7\textwidth]{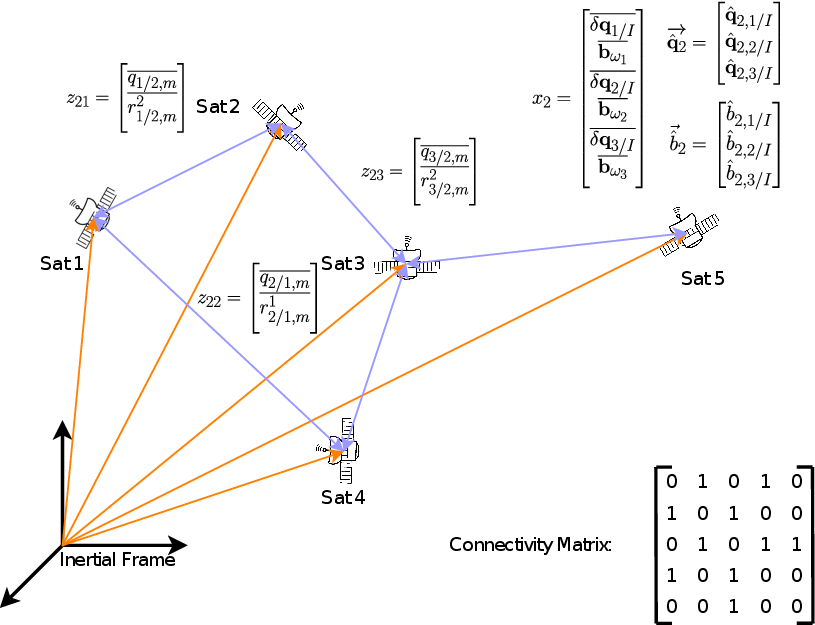}
\caption{Configuration of a fleet of $l=5$ satellites.
The states of satellite 2 are shown, where $x_2$ is the filter state, and $\protect\overrightarrow{\hat{\mathbf{q}_2}}$ and $\protect\overrightarrow{\hat b_2}$ are the respective estimated pose and bias of satellite 2 and its neighbours.}
\label{fig:problem_statement}
\end{figure*}

The relative pose measurements can be taken by a visual sensor \cite{visnav,visnav2,kreiss2021openpifpaf,Kreiss_2019_CVPR}. 
The absolute pose measurements can be taken by a combination of star trackers and a Global Positioning System (GPS). 
The angular velocity and linear acceleration measurement can be taken via an IMU. 
Furthermore, we assume that each satellite has communication links to its neighbouring satellites for data exchange. 
Each satellite therefore only keeps track of its own state, and the state of its neighbours.

\subsection{Derivation of the distributed DQ-MEKF (DDQ-MEKF)}
\label{sec:derivation}

In this setting, consider the set of satellites\footnote{We will use `nodes' and `satellites' interchangeably from this point onwards} ${V}$ with $|{V}| = l$. 
Each satellite $i\in V$ will only estimate its own pose and the pose of its neighbours, i.e., it can only estimate the poses of satellites  $j\in{V}_i := N_i \cup i$. 
Furthermore, each node $i$ shall compute the pose, the pose error, as well as the velocity measurement bias and velocity of the $i$th satellite. 
We can thus define the state and process noise vector of each satellite.

Let $i \in {V}$.
We write ${V}_i = \{j_{i, 1},\dots,j_{i, k_i +1}\} \subset {V}$, where $j_{i,k} \in {V}_i$ is the $k$th neighbour of $i$ -- note that due to the definition of $V_i$, satellite $i$ is also treated as a neighbour of $i$.
Since $i$ knows the node labels of its neighbours, we can assume without loss of generality that ${V}_i$ is sorted in ascending order of those labels. 
Then, we denote for each $i=1,\dots, l$ the state variables $x_i \in \mathbb{R}^{12(k_i+1)}$ and $w_i  \in \mathbb{R}^{12(k_i+1)} $ given by,
\begin{align}
  \begin{split}
    x_i &= \begin{bmatrix}
      \overline{\dq{\delta q}}_{j_{i, 1}/I}^T & 
      \overline{\dq{b_\omega}}_{j_{i, 1}}^T & 
      \cdots & 
      \overline{\dq{\delta q}}_{j_{i, k_i+1}/I}^T & 
      \overline{\dq{b_\omega}}_{j_{i, k_i+1}}^T
    \end{bmatrix}^T
     \\
    w_i &= \begin{bmatrix}
      \overline{\dq{\eta_{\omega}}}_{j_{i,1}}^T & 
      \overline{\dq{\eta_{b_\omega}}}_{j_{i, 1}}^T & 
      \cdots & 
      \overline{\dq{\eta_{\omega}}}_{j_{i, k_i+1}}^T & 
      \overline{\dq{\eta_{b_\omega}}}_{j_{i, k_i+1}}^T 
    \end{bmatrix}^T,
  \end{split}\label{eq:distributed_kf_state_def}
\end{align}
where we
simplify the notation via $\overline{\dq{\delta q}}_{j_{i, k}/I} := \overline{\dq{\delta q}}_{B_{j_{i, k}/I}}$.
Following the centralized filter described in 
\if\arxivyes1
Appendix~\ref{sec:appendix-1},
\else
Appendix B of the supplementary material~\cite{arxivversion},
\fi
  $ \overline{\dq{\delta q}}_{j_{i, k}/I}$ is the error dual quaternion representing the pose error relative to a reference of the $k$th neighbour of $i$, $\overline{\dq{b_\omega}}_{j_{i, k}}$ is the dual bias of the velocity, and $w_i$ is a vector of process noises.
 
In contrast to the non-cooperative case in 
\if\arxivyes1
Appendix~\ref{sec:dq-mekf},
\else
Appendix B of the supplementary material~\cite{arxivversion},
\fi
 the definition of the state and noise in Eqn~\eqref{eq:distributed_kf_state_def} is  a stacking of the states of satellite $i$ and its neighbours. 
The top state is always the state of the satellite with the smallest nodal label. 
We further define the stacked pose estimate vector $\overrightarrow{\hat{\mathbf{q}}_i}$ and the stacked dual bias vector $\overrightarrow{\hat b_i}$ for $i=1,\dots, l$ as,
\begin{subequations}\label{eq:distributed_kf_pose_bias_def}
  \begin{align}
    \overrightarrow{\hat{\mathbf{q}}_i} &= 
    \begin{bmatrix}
      \dq{\hat q}_{i, j_{i, 1}/I}^T &
      \cdots &
      \dq{\hat q}_{i, j_{i, k_i+1}/I}^T 
    \end{bmatrix}^T
    \in \mathbb{R}^{8(k_i+1)} \\
    \overrightarrow{\hat b_i} &= 
    \begin{bmatrix}
      \overline{\dq{\hat b_{\omega}}}_{i, j_{i, 1}}^T &
      \cdots &
      \overline{\dq{ \hat b_{\omega}}}_{i, j_{i, 1}}^T
    \end{bmatrix}^T
    \in \mathbb{R}^{6(k_i+1)}.
  \end{align}  
\end{subequations}
Again, to simplify notation we denote $\dq{\hat q}_{i, k/I} := \dq{\hat q}_{i, B_k/I}$ as the estimate that satellite $i$ computes of the pose of neighbouring satellite $k$.

\subsubsection{Time Update}\label{sec:time_update_dist}

At the beginning of each iteration, each satellite $i \in V$ takes a measurement $\dq{\omega}_{B_i/I, m}^{B_i}$ (where the subscript $m$ denotes `measurement') of its own dual velocity in its body frame $B_i$, and propagates this information to its neighbours. 
To reduce notational burden of notation, we write $\dq{\omega}_{j_{i,k}/I,m}:=\dq{\omega}_{B_{j_{i, k}}/I, m}^{B_{j_{i, k}}}$ to denote the measured dual velocity of the $k$th neighbour of satellite $i$ in the body frame of the $k$th neighbour of satellite $i$.
Hence, each satellite $i\in V$ can construct the stacked dual velocity measurement as,
\begin{equation}
    \label{eq:stacked_vel_measurement_vec}
    \overrightarrow{\dq{\omega}}_{i, m} = 
    \begin{bmatrix}
      \overline{\dq{\omega}_{j_{i,1}/I,m}}^T&  
      \cdots &
      \overline{\dq{\omega}_{j_{i,k_i+1}/I,m}}^T
    \end{bmatrix}^T
    \in \mathbb{R}^{6(k_i+1)}.
\end{equation}
By removing the dual bias from the dual velocity, each $i\in V$ can compute the estimated stacked dual velocity via,
\begin{equation}
    \overrightarrow{ \hat{\dq{\omega}}}_{i} = \overrightarrow{\dq{\omega}}_{i, m} - \overrightarrow{\hat b_i} = 
    \begin{bmatrix}
      \overline{\hat{\dq{\omega}}_{j_{i,1}/I,m}}^T&  
      \cdots &
      \overline{\hat{\dq{\omega}}_{j_{i,k_i+1}/I,m}}^T
    \end{bmatrix}^T
\label{eq:stacked_vel_est_vec}
\end{equation}
This yields the set of differential equations for the estimate,
\begin{equation}\label{eq:distributed_kf_propagation}
  \dfrac{d}{dt}
 \begin{bmatrix}
\dq{\dot{\hat{q}}}_{i, {j_{i, 1}/I}}\\
\vdots\\
\dq{\dot{\hat{q}}}_{i, {j_{i, k_i+1}}/I}
\end{bmatrix} \approx
\begin{bmatrix}
\frac{1}{2}\dq{\hat q}_{i, {j_{i, 1}}/I}\hat{\dq{\omega}}_{j_{i,1}/I,m}\\
\vdots\\
\frac{1}{2}\dq{\hat q}_{i, {j_{i, k_i+1}}/I}\hat{\dq{\omega}}_{j_{i,k_i+1}/I,m}
\end{bmatrix} = \frac{1}{2}
\begin{bmatrix}
\lqm{\dq{\hat q}_{i, {j_{i, 1}}/I}} & & 0\\
& \ddots & \\
0 & & \lqm{\dq{\hat q}_{i, {j_{i, k_i+1}}/I}}
\end{bmatrix}
\begin{bmatrix}
\hat{\dq{\omega}}_{j_{i,1}/I,m}\\
\vdots\\
\hat{\dq{\omega}}_{j_{i,k_i+1}/I,m}
\end{bmatrix},
\end{equation}
from which the time update of the state is computed. 

Next, the covariance of the state is propagated.
Each node $i\in V$ must compute the Jacobians of the dynamics $F_{i,m}$ and $G_{i,m}$ for each of its neighbours  $m\in{V}_i$ using Eq.~\eqref{eq:F_and_G_def}.
Explicitly, these are, 
\begin{subequations}\label{eq:distributed_Fik_Gik}
\begin{align}
    F_{i, m} &= \begin{bmatrix}-\overline{\hat{\dq{\omega}}_{{i, j_{i, k}}/I}}^{\times} & -\frac{1}{2}\mathbf{I}_{6\times 6} \\ \mathbf{0}_{6\times 6} & \mathbf{0}_{6\times 6}\end{bmatrix}  \text{for }m=1,\dots k_i+1,\\
    G_{i, m} &= \begin{bmatrix}-\frac{1}{2}\mathbf{I}_{6 \times 6} & \mathbf{0}_{6 \times 6} \\ \mathbf{0}_{6 \times 6} & \mathbf{I}_{6\times 6}\end{bmatrix} \quad \text{for }m=1,\dots k_i+1.
    \end{align}
\end{subequations}
The complete Jacobians of the dynamics and measurement models $F_i$ and $G_i$ of satellite $i\in V$ are computed by concatenating the aforementioned $F_{i,k}$ and $G_{i,k}$ matrices block-diagonally, yielding
\begin{subequations}\label{eq:distributed_F_G}
  \begin{align}
    F_i &= \mathrm{blkdiag}\{F_{i,m}\}_{m=1}^{k_i+1}\in \mathbb{R}^{12(k_i+1)\times 12(k_i+1)} \\
    G_i &= \mathrm{blkdiag}\{G_{i,m}\}_{m=1}^{k_i+1}\in \mathbb{R}^{12(k_i+1)\times 12(k_i+1)}.
  \end{align}
\end{subequations}
The stacked $Q_i$ matrix therefore becomes,
\begin{align}
  Q_i &= \mathrm{blkdiag}\left\{\tilde{Q}_{i,m}\right\}_{m=1}^{k_i+1}\\
  \tilde{Q}_{i,m} &= \mathrm{blkdiag}\left\{\overline{Q}_{\dq{\omega}_{j_{i, m}}}, \overline{Q}_{\dq{{b_\omega}}_{j_{i, m}}}\right\}.
\end{align}
Next, each node $i\in V$ can propagage the covariance using Eq.~\eqref{eq:Ricatti}
As stated above, the only communication needed between satellites is the exchange of the dual velocity measurement $\dq{\omega}_{B_i/I, m}^{B_i}$.
In the case where this information is not accessible (i.e., no gyro or linear acceleration measurement), no communication is necessary.

\subsubsection{Measurement Update}
We must fuse both the relative and absolute pose measurements in order to derive the measurement update step. 
As explained in \S\ref{sec:problem_statement}, each satellite $i\in V$ has access to an absolute pose measurement $\dq{q}_{B_i/I, m}$, and $k_i$ relative attitude and relative position measurements, $q_{B_k/B_i, m}$ and $r_{B_k/B_i, m}^{B_i}$, in its own body frame.
We can thus denote the measurement vector $z_i$ of satellite $i$ as,
\begin{align}
    z_i &= 
    \begin{bmatrix}
        \overline{q_{B_{n_{i,1}}/B_i, m}}^T&
        \overline{r_{B_{n_{i,1}}/B_i, m}^{B_i}}^T&
        \cdots&
        \overline{\dq{\hat q}_{i,i}^*\dq{q}_{B_i/I, m}}^T&
        \cdots&
        \overline{q_{B_{n_{i,k_i}}/B_i, m}}^T&
        \overline{r_{B_{n_{i,k_i}}/B_i, m}^{B_i}}^T&
    \end{bmatrix}^T\\ &:= 
        \begin{bmatrix}
             \overline{h_{q_{n_{i,1}/i}}(x_i)}^T&
             \overline{h_{r_{n_{i,1}/i}^i}(x_i)}^T&
             \cdots&
             \overline{h_{abs}(x_i)}^T&
             \cdots&
             \overline{h_{q_{n_{i,k_i}/i}}(x_i)}^T&
             \overline{h_{r_{n_{i,k_i}/i}^i}(x_i)}^T&
    \end{bmatrix}^T + 
    v_i.\label{eq:measurement_vec_distr}
\end{align}
Note that $z_i\in \mathbb{R}^{6(k_i+1)}$.
In Equation \eqref{eq:measurement_vec_distr}, $h_{q_{n_{i,k}/i}}$ is the function which outputs the relative measurement from satellite $i$ and its neighbour with index $n_{i,k}$ to the state $x_i$ of satellite $i$. 
In a similar manner, $h_{r_{n_{i,k}/i}^i}$ gives the relative position. 
Only one entry is different, namely $\overline{\dq{\hat q}_{i,i}^*\dq{q}_{B_i/I, m}}$, which is the satellite's measurement of its own state.
$v_i$ is the measurement noise acting on all measurements.

The Jacobian $H$ of the measurement function given in Eq.~\eqref{eq:measurement_vec_distr} can be viewed as the (linearized around zero error) mapping from the state to the measurement.
In the case of Eq.~\eqref{eq:measurement_vec_distr}, the Jacobian consists of three submatrices: the mapping of the absolute pose measurement to the state, the mapping from the relative attitude measurement to the state, and finally the mapping from the relative position measurement to the state. 
For sake of clarity, we derive the Jacobians of the individual measurement functions $h$ in Appendix~\ref{sec:jacob-meas-models}, and proceed to the construction of the full Jacobian of the measurement~in Eq.~\eqref{eq:measurement_vec_distr}.

\textit{Full Jacobian of the Measurement:}\\
We begin as follows.
We define the following four matrices as,
\begin{align}
\overline{H^i_{i, k}} &:= \begin{bmatrix}
                \frac{1}{2}\overline{H_{q_{i/k}, \delta q_{i/I}}} & \mathbf{0}_{3\times 3}\\
                \frac{1}{4}\overline{H_{r_{k/i}^i, \delta q_{i/I}}} & \frac{1}{4}\overline{H_{r_{k/i}^i, \delta p_{i/I}}} 
            \end{bmatrix},\quad
\overline{H^k_{i, k}} := \begin{bmatrix}
                \frac{1}{2}\overline{H_{q_{i/k}, \delta q_{k/I}}} & \mathbf{0}_{3\times 3}\\
                \frac{1}{4}\overline{H_{r_{k/i}^i, \delta q_{k/I}}} & \frac{1}{4}\overline{H_{r_{k/i}^i, \delta p_{k/I}}}
            \end{bmatrix}\\
H^i_{i, k} &:= \begin{bmatrix}
                \frac{1}{2}H_{q_{i/k}, \delta q_{i/I}} & \mathbf{0}_{4\times 4}\\
               \frac{1}{4} H_{r_{k/i}^i, \delta q_{i/I}} & \frac{1}{4}H_{r_{k/i}^i, \delta p_{i/I}} 
            \end{bmatrix}, \quad
H^k_{i, k} := \begin{bmatrix}
                \frac{1}{2}H_{q_{i/k}, \delta q_{k/I}} & \mathbf{0}_{4\times 4}\\
                \frac{1}{4}H_{r_{k/i}^i, \delta q_{k/I}} & \frac{1}{4}H_{r_{k/i}^i, \delta p_{k/I}}
            \end{bmatrix},
\end{align}
where the submatrices $H_a\in\mathbb{R}^{4\times 4},\overline{H}_a\in\mathbb{R}^{3\times 3}$ corresponding to the Jacobians of the various measurement models of attitude and position are derived in Appendix~\ref{sec:jacob-meas-models}.
The matrix $H^i_{i, k}$ can be understood as the part of the measurement Jacobian which transforms the state of satellite $i$ to the relative measurement that satellite $i$ takes of satellite $k$. 
On the other hand, $H^k_{i, k}$  is the part of the measurement Jacobian which transforms the state of satellite $k$, estimated from satellite $i$ to the relative measurement that satellite $i$ takes of satellite~$k$.
Here, we define $\overline{H_a}\in \mathbb{R}^{3\times 3}$ as the matrix $H_a \in \mathbb{R}^{4\times 4}$ with its first row and column (corresponding to the scalar part of the quaternion) deleted.
Thus, $\overline{H_a}$ is congruent with the reduced state space defined in Eq.~\eqref{eq:distributed_kf_state_def}. 

Finally, let $H_i := H_{i,i}$ and $\overline{H_i} := \overline{H_{i,i}}$.
The stacked measurement Jacobian for each satellite $i$ is given by,
\begin{align}
        H_i &= \begin{bmatrix}
            H^{n_{i, 1}}_{i, n_{i, 1}} &  \mathbf{0}_{8 \times 6} & \dots & H^i_{i, n_{i, 1}} & \mathbf{0}_{8 \times 6} & \dots & \mathbf{0}_{8 \times 8} & \mathbf{0}_{8 \times 6}\\
            \vdots & \vdots & \ddots & \vdots & \vdots & \ddots & \vdots & \vdots\\
            \mathbf{0}_{8 \times 8} & \mathbf{0}_{8\times 6} & \dots & \mathbf{I}_{8\times 8} & \mathbf{0}_{8 \times 6} & \dots & \mathbf{0}_{8 \times 8} & \mathbf{0}_{8 \times 6}\\
            \vdots & \vdots & \ddots & \vdots & \vdots & \ddots & \vdots & \vdots\\
             \mathbf{0}_{8 \times 8} & \mathbf{0}_{8 \times 6} & \dots & H^i_{i, n_{i, k_i}} & \mathbf{0}_{8 \times 6} & \dots & H^{k_i}_{i, n_{i, k_i}} &  \mathbf{0}_{8 \times 6} 
            \end{bmatrix}\in \mathbb{R}^{8(k_i+1)\times 14(k_i+1)}\\
        \overline{H_i} &= \begin{bmatrix}
            \overline{H^{n_{i, 1}}_{i, n_{i, 1}}} &  \mathbf{0}_{6 \times 6} & \dots & \overline{H^i_{i, n_{i, 1}}} & \mathbf{0}_{6 \times 6} & \dots & \mathbf{0}_{6 \times 6} & \mathbf{0}_{6 \times 6}\\
            \vdots & \vdots & \ddots & \vdots & \vdots & \ddots & \vdots & \vdots\\
            \mathbf{0}_{6 \times 6} & \mathbf{0}_{6\times 6} & \dots & \mathbf{I}_{6\times 6} & \mathbf{0}_{6 \times 6} & \dots & \mathbf{0}_{6 \times 6} & \mathbf{0}_{6 \times 6}\\
            \vdots & \vdots & \ddots & \vdots & \vdots & \ddots & \vdots & \vdots\\
             \mathbf{0}_{6 \times 6} & \mathbf{0}_{6\times 6} & \dots & \overline{H^i_{i, n_{i, k_i}}} & \mathbf{0}_{6 \times 6} & \dots & \overline{H^{k_i}_{i, n_{i, k_i}}} &  \mathbf{0}_{6 \times 6} 
            \end{bmatrix}\in \mathbb{R}^{6(k_i+1)\times 12(k_i+1)},\label{eq:distributed_H}
\end{align}
where $\overline{H^{n_{i, j}}_{i, n_{i, j}}}$ is always in the $j$-th column and $j$-th row if $n_{i, j} < i$, and in the $j+1$-th column and $j+1$-th row otherwise. 
The matrix $\overline{H^i_{i, n_{i, j}}}$ is always in the $i$-th column and $j$-th row if $n_{i, j} < i$, and in the $i$-th column and $j+1$-th row otherwise. 
The matrix $\mathbf{I}_{6 \times 6}$ is always in the $j$-th column and row where $n_{i, j} = i$.

To complete the measurement update, we must compute the measurement covariance matrix as well as the transformed data vector required for the information filter.
These are given by,

\begin{align}
    \overline{U_i} &= \overline{H_i^\intercal}\overline{R_i}^{-1}\overline{H_i}\label{eq:meas_cov_distr1} \\
  u_i &= \overline{H_i^\intercal}\overline{R_i}^{-1}z_i \label{eq:meas_cov_distr2}\\
  U_i &= H_i^\intercal R_i^{-1}H_i,\label{eq:meas_cov_distr3}
\end{align}
where,
\begin{align}
  \overline{R}_i &:= \Cov[v_i] = \mathrm{blkdiag}\left\{ \overline{R_{n_{i,m},i}} \right\}_{m=1}^{k_i+1}\\
  {R}_i &:= \Cov[v_i] = \mathrm{blkdiag}\left\{ {R_{n_{i,m},i}} \right\}_{m=1}^{k_i+1}.
\end{align}
Here, $\overline{R_{n_{i,m}, i}} \in \mathbb{R}^{6\times 6}$ is the covariance of the measurement of the pose of satellite $m$ as taken by satellite $i$ in Equation~\eqref{eq:measurement_vec_distr}. 
The matrix $R_{n_{i,m}, i}$ is the `extension' of $\overline{R_{n_{i,m}, i}}$ to $\mathbb{R}^{8 \times 8}$, where a row and column are added corresponding to each scalar entry of the dual quaternion.
We add entries to the diagonal (say, 1) in these new rows and columns such that the matrix is not singular, however, what is set in those two added entries is of no importance as those entries will not be propagated. 
Further note that $\overline{R_i}\in \mathbb{R}^{6(k_i+1)\times 6(k_i+1)}$ and $ {R_i}\in \mathbb{R}^{8(k_i+1)\times 8(k_i+1)}$.

Next, the Kalman state update is computed by each node $i\in V$ via,
\begin{align}
    \Delta x_i &= 
    \begin{bmatrix}
    \Delta \overline{\dq{\delta \hat{q}}_{j_{i,1}/I}}^T &
    \Delta \overline{\dq{\delta \hat{b}}_{\omega_{j_{i,1}}}}^T&
    \cdots &
    \Delta \overline{\dq{\delta \hat{q}}_{j_{i,k_i+1}/I}}^T&
    \Delta \overline{\dq{\delta \hat{b}}_{\omega_{j_{i,k_i+1}}}}^T&
    \end{bmatrix} 
                 = M_i\left(\overline{u_i} - \overline{\left(U_i \mathbf{1}\right)}\right),\label{eq:dqekf-update}
\end{align}
where $M_i$ is given by Eq.~\eqref{eq:M_matrix} as,
\begin{equation}\label{eq:M_matrix_distributed}
    M_i = \left(P_i^{-1} + \overline{U_i}\right)^{-1}.
\end{equation}
Finally, the Kalman state update is used to compute the time update.
Similarly as to the non-cooperative case in 
\if\arxivyes1
Appendix~\ref{sec:dq-mekf},
\else
the supplementary material~\cite{arxivversion},
\fi
 each node $i\in V$ updates each dual quaternion entry in $\Delta x_i$ for $k=1,\dots k_i+1$ as, 
\begin{align}
  \label{eq:extend_to_R8_dist}
  \begin{split}
    \Delta \delta \hat q_{j_{i,k}, r} &= \left(\sqrt{1-||\Delta\overline{\delta \hat{q}_{j_{i,k}, r}}||^2}, \Delta\overline{\delta \hat{q}_{j_{i,k}, r}}\right)  \\
\Delta \delta \hat q_{j_{i,k}, d} &= \left(\frac{-\Delta\overline{\delta \hat{q}_{j_{i,k}, r}}^T\Delta\overline{\delta \hat{q}_{j_{i,k}, d}}}{\sqrt{1-||\Delta\overline{\delta \hat{q}_{j_{i,k}, r}}||^2}}, \Delta\overline{\delta \hat{q}_{j_{i,k}, d}}\right)
  \end{split}
\end{align}
and again each node $i\in V$ concludes  the measurement update by computing the estimate for $k=1,\dots k_i+1$  as,
\begin{align}
  \label{eq:dist_meas_update_estimated}
  \begin{split}
    \dq{\hat q}^+_{i, k/I} &= \dq{\hat q}^-_{i, k/I} \Delta \delta \dq{\hat q}_{j_{i,k/I}} \\
    \overline{ \dq{\hat b}^+_{\omega_{i, k}}} &= \overline{ \dq{\hat b}^-_{\omega_{i, k}}} + \Delta \overline{\delta \dq{\hat b}_{\omega_{i, k}}}.
  \end{split}
\end{align}
The DDQ-MEKF is summarized in Algorithm~\ref{alg:ddq-mekf}.

\begin{algorithm}
  \caption{Distributed Dual Quaternion Kalman Filter}
  \label{alg:ddq-mekf}
  \begin{algorithmic}
    \Require $\forall i \in J: Q_i,R_i,\overrightarrow{\hat q_i}(t_0):=q_0,\hat P_i(t_0) := P_0$
    \While{$t_k \geq 0$}
    \For{$i\in J$}
    \State Take measurement $\dq{\omega}^{B_i}_{B_i/I, m}$
    \State Propagate measurement to neighbours $j\in N_i$
    \State Compute dual velocity estimate $\overrightarrow{\hat \omega_i} = \overrightarrow{\omega}_{i, m} - \overrightarrow{\hat b_i}$
    \State Propagate each estimate $\overrightarrow{\hat q^+_i} \rightarrow \overrightarrow{\hat q^-_i}$ via Eqn.~\eqref{eq:distributed_kf_propagation}
    \State Compute $F$, $G$; propagate $P$ via Eqns.~\eqref{eq:distributed_F_G}, \eqref{eq:Ricatti}
    \State Measure own and relative poses $z_{i}$ via Eqn.~\eqref{eq:measurement_vec_distr}
    \State Compute $\overline{H_i}$ and $H_i$ matrix via Eqn.~\eqref{eq:distributed_H}
    \State Compute $\overline{U_i}, U_i,$ and $\overline{u_i}$ via Eqn.~\eqref{eq:meas_cov_distr1},~\eqref{eq:meas_cov_distr2} \&~\eqref{eq:meas_cov_distr3}
    \State Compute $M_i$ via Eqn.~\eqref{eq:M_matrix_distributed}
    \State Compute $\Delta x_i$ via Eqn.~\eqref{eq:dqekf-update}
    \State Extend  Kalman state update via Eqn.~\eqref{eq:extend_to_R8_dist}
    \State Update estimate and bias via Eqn.~\eqref{eq:dist_meas_update_estimated}  
    \EndFor
    \EndWhile
  \end{algorithmic}
\end{algorithm}

\subsection{Consensus}
In Algorithm~\ref{alg:ddq-mekf}, the satellites communicate with each other for velocity measurements only.
In this subsection, we describe a consensus algorithm that allows each satellite to update and improve its state estimate by incorporating the state estimates of its neighbours.

We propose two consensus algorithms, denoted \emph{soft} and \emph{hard}. 
The goal of soft consensus is to average the state estimates of neighbouring satellites without any knowledge of the error covariance.
Hard consensus, building upon~\cite{1583238}, fuses covariance data exchanged between neighbours. 

\subsubsection{Description of the Soft Consensus Term}\label{sec:description_soft_cons}
Before formally stating the algorithm, we give an example which illustrates the soft consensus step.
Consider satellite 4 in Fig.~\ref{fig:problem_statement}. 
Satellite 4 has two neighbours, namely satellites 1 and 3, and thus has no notion of any other satellites in the fleet. 
After the measurement update of each satellite, the available estimates are listed in Table~\ref{tab:estimates}.
\begin{table}
  \centering
  \begin{tabular}{ c c c }
Satellite 1 & Satellite 3 & Satellite 4\\\hline\\
$\overrightarrow{\hat{\mathbf{q}}^+_1} = \begin{bmatrix}\dq{\hat q}^+_{1, 1/I} \\ \dq{\hat q}^+_{1, 2/I} \\ \dq{\hat q}^+_{1, 4/I} \end{bmatrix}$ &
$\overrightarrow{\hat{\mathbf{q}}^+_3} = \begin{bmatrix}\dq{\hat q}^+_{3, 2/I} \\ \dq{\hat q}^+_{3, 3/I} \\ \dq{\hat q}^+_{3, 4/I} \\ \dq{\hat q}^+_{3, 5/I} \end{bmatrix}$ &
$\overrightarrow{\hat{\mathbf{q}}^+_4} = \begin{bmatrix}\dq{\hat q}^+_{4, 1/I} \\ \dq{\hat q}^+_{4, 3/I} \\ \dq{\hat q}^+_{4, 4/I} \end{bmatrix} $
\end{tabular}
\caption{Example of available estimates to each satellite, following Fig.~\ref{fig:problem_statement}.}\label{tab:estimates}
\end{table}
Hence, satellite 4 can use certain measurements from its neighbour to improve its own state estimates. 
From satellite 1, it can use $\dq{\hat q}^+_{1, 1}$ and $\dq{\hat q}^+_{1, 4}$; similarly, from satellite 3 it can use $\dq{\hat q}^+_{3, 3}$ and $\dq{\hat q}^+_{3, 4}$. 
The other estimates from the neighbouring satellites are not used, as those are estimates of states of satellites not in the neighbourhood of satellite 4. 
The same notion pertains to the estimated biases.
 
The idea is now to average the estimated states of satellite 4, namely $\overrightarrow{\hat q^+_4}$, with the estimated states of satellite 1 and 3 which they have in common. 
In other words we would like to implement the consensus term in Eq.~\eqref{eq:olfati_state_update},
\begin{equation}\label{eq:distributed_soft_cons_def}
    \phi_i = \mu_i \sum_{j\in N_i} \left(x^-_j - x^-_i\right),
\end{equation}
adapted to the dual quaternion setting. 
As our state estimate is the \emph{error} of the quaternion, it is nonsensical to use Eq.~\eqref{eq:distributed_soft_cons_def}.
Instead, the Kalman state update is computed first, then the pose estimate is updated, and then a consensus step is performed. 
Furthermore, computing the arithmetic difference of the pose as is done in Eq.~\eqref{eq:distributed_soft_cons_def} must be adapted to the dual quaternion setting; hence, we propose a consensus step which utilizes a more natural notion of  averaging dual quaternions.

\textit{Soft Consensus on the Position and Velocity:}\\
After performing the measurement update, each satellite $i\in V$ has an estimate of its own pose and bias and the pose and bias of its neighbours:
\begin{align}
\label{eq:distributed_kf_pose_bias_def}
    \overrightarrow{\hat{\mathbf{q}}^+_i} &= 
  \begin{bmatrix}
    \dq{\hat q}^{+,T}_{i, j_{i, 1}/I} & 
    \cdots & 
    \dq{\hat q}^{+,T}_{i, j_{i, k_i+1}/I}
  \end{bmatrix}^T \in \mathbb{R}^{8(k_i+1)}\\
    \overrightarrow{\hat b^+_i} &= 
  \begin{bmatrix}
  \overline{\dq{\hat b}^{+,T}_{\omega}}_{i, j_{i, 1}} 
   & 
   \cdots 
   & 
   \overline{\dq{ \hat b}^{+,T}_{\omega}}_{i, j_{i, 1}}
 \end{bmatrix}^T \in \mathbb{R}^{6(k_i+1)}.
\end{align}
Each satellite $i\in V$ first needs to propagate this data to its neighbours $k\in N_i$:
\begin{equation}
    \text{msg}_{i,k} = \left( \overrightarrow{\hat{\mathbf{q}}_i^+}, \overrightarrow{\hat b_i^+}\right).
\end{equation}
We now focus on the issue of the correct arithmetic notion of the `difference' between two dual quaterions.
First, the dual part of the dual quaternion is transformed to the position using $r_{B/I}^B = \begin{bmatrix} 0 & \overline{r_{B/I}^B}^T\end{bmatrix}^T = 2p_{B/I}q_{B/I}$, yielding the position estimate for each $i \in V$, 
\begin{equation}\label{eq:distributed_kf_position_def}
\overline{\hat r^+_i} = 
\begin{bmatrix}
\overline{\hat r^{+,T}_{i, j_{i, 1}}} & 
\cdots &  
\overline{\hat r^{+,T}_{i, j_{i, k_i+1}}} 
\end{bmatrix}^T \in \mathbb{R}^{3(k_i+1)}.
\end{equation}
Since each satellite is estimating its own state and that of its neighbours, the various state vectors of each satellite are not compatible with one another.
To rectify this,  we extend each satellite's estimated position vector to $\mathbb{R}^{3|V|} = \mathbb{R}^{3l}$ by padding the entries of untracked satellites with zeros.
Then, we copy this position vector for each neighbour $k\in N_i$ of satellite $i$, and delete position estimates that satellite $i$ and $k$ do not have in common.

Formally, we define the subset $\Lambda_{i, k} := \left({V}_i \cap {V}_k\right) \subset {V}$, in other words the set of nodes which satellite $i$ and $k$ have in common. 
For example, satellite 2 and 3 in the running example of Fig.~\ref{fig:problem_statement} have $\Lambda_{2, 3} = \{2, 3\}$.  
We write $\Lambda_{i, k} = \{\lambda_{ik, 1},\dots,\lambda_{ik, |\Lambda_{i, k}|}\} \subset {V}$, i.e.~$\Lambda_{i, k}$ is sorted in ascending order, as again, we assume that satellites know their own label in $V$.
The full position estimate vector of satellite $i\in V$ with respect to its neighbour $k\in N_i$ is thus defined as,
\begin{align}\label{eq:distributed_kf_position_def_expanded}
    \overrightarrow{\hat r^+}_{\text{full}, i,k} &= 
  \begin{bmatrix} 
    \overline{r_1}^T & \cdots & \overline{r_l}^T 
  \end{bmatrix}^T \in \mathbb{R}^{3l} \\ 
\text{ where } \overline{r_j} &= \begin{cases}
    \overline{\hat r^+_{i, j}},      & \text{if } j \in \Lambda_{i,k}\\
    \mathbf{0}_{3\times 1},              & \text{otherwise},
\end{cases}
\end{align}
and where $r^+_{i, j}$ is the estimated position of satellite $j$ by satellite $i$. 

We can now easily subtract position vectors between neighbouring satellites to compute their difference: 
\begin{equation}\label{eq:distributed_soft_cons_pos}
    \phi_{\text{full}, i, r} = \mu_{i, r} \sum_{j\in N_i} \left(\overrightarrow{\hat r^+}_{\text{full}, j,i} - \overrightarrow{\hat r^+}_{\text{full}, i,j}\right) \in \mathbb{R}^{3l}.
\end{equation}
Now, $\phi_{\text{full}, i, r}$ needs to be reduced to remove the padded zeros corresponding to positions that satellite $i$ does not estimate.
This is achieved by removing each entry corresponding to the indices of $\overline{r_j}$ in Eq.~\eqref{eq:distributed_kf_position_def_expanded} which satisfies $ j\notin {V}_i$, yielding $\phi_{i, r} \in \mathbb{R}^{3(k_i+1)}$.
This concludes the portion of the soft consensus step on the position part of the dual quaternion.

For velocity consensus, as one can see from Eq.~\eqref{eq:stacked_vel_est_vec} it is sufficient to perform consensus on the bias.
Therefore, achieving consensus on the bias automatically achieves consensus on the velocity. 
Consensus on the bias is performed using the same procedure as on position. 
We define the full bias vector of satellite $i\in V$ with respect to its neighbour $k\in N_i$ as,
\begin{align}\label{eq:distributed_kf_bias_def_expanded}
    \overrightarrow{\hat b^+}_{\text{full}, i,k} &= 
       \begin{bmatrix} 
           b_1^T & \cdots & b_l^T 
       \end{bmatrix}^T \in \mathbb{R}^{6l} \\\text{ where } b_j &= \begin{cases}
   \overline{\dq{\hat b}^+_{\omega}}_{i, j},      & \text{if } j \in \Lambda_{i,k}\\
    \mathbf{0}_{6 \times 1},              & \text{otherwise},
\end{cases}
\end{align}
and where $\overline{\dq{\hat b}^+_{\omega}}_{i, j}$ is the estimated bias of satellite $j$ by satellite $i$. 
The consensus step on the bias is therefore,
\begin{equation}\label{eq:distributed_soft_cons_bias}
    \phi_{\text{full}, i, b} = \mu_{i, b} \sum_{j\in N_i} \left(\overrightarrow{\hat b^+}_{\text{full}, j,i} - \overrightarrow{\hat b^+}_{\text{full}, i,j}\right) \in \mathbb{R}^{6l}.
\end{equation}
The resulting vector $\phi_{\text{full}, i, b}$ is then reduced by removing the extraneous entries corresponding to the padded zeros, as was done on $\phi_{\text{full}, i, r}$, yielding $\phi_{i, b} \in \mathbb{R}^{6(k_i+1)}$.

\textit{Soft Consensus on the Attitude}\\
Next, we define consensus on the attitude quaternion by exploiting a notion of `averaging' a quaternion.
First, for each $i\in V$, we define the attitude vector as,
\begin{equation}\label{eq:distributed_kf_att_def}
\overline{\hat q^+_i} = 
\begin{bmatrix}
  \hat q^{+,T}_{i, j_{i, 1}} & 
  \cdots & 
  \hat q^{+,T}_{i, j_{i, k_i+1}} 
\end{bmatrix}^T \in \mathbb{R}^{4(k_i+1)}.
\end{equation}
Next, as for the position vector case, we  pad the attitude vector with additional elements for the parts of the attitude corresponding to satellites which $i$ doesn't track, and create copies of the vector for each neighbour $k\in N_i$ which only contain estimates in common:
\begin{align}\label{eq:distributed_kf_attitude_def_expanded}
    \overrightarrow{\hat q^+}_{\text{full}, i,k} &= 
  \begin{bmatrix} 
    q^{i,k,T}_1 & 
    \cdots & 
    q^{i,k,T}_l 
  \end{bmatrix}^T \in \mathbb{R}^{4l} \\ 
  \text{ where } q^{i,k}_j &= \begin{cases}
   \hat q^+_{i, j},      & \text{if } j \in \Lambda_{i,k}\\
    \begin{bmatrix}1 & 0 & 0 & 0\end{bmatrix}^T,              & \text{otherwise},
\end{cases}
\end{align}
and where $q^+_{i, j}$ is the estimated attitude of satellite $j$ by satellite $i$. 

The correct notion of the `difference' between quaternions $q_1$ and $q_2$ can be expressed as $q_1^*q_2$, and the correct notion of `summation' over all of the quantites is given by quaternion multiplication.
Hence, the correct notion of the consensus step for the attitude quaternions is given by,
\begin{align}
    \theta_{\text{full}, i, q} &= \prod_{j\in N_i} \left(\overrightarrow{\hat q^+}_{\text{full}, i,j}\right)^*\left(\overrightarrow{\hat q^+}_{\text{full}, j,i}\right)\label{eq:prep_distributed_soft_cons_attitude}\\
                  &= \prod_{j\in N_i} \begin{bmatrix} \left(q^{i,j}_1\right)^*\left(q^{j,i}_1\right) \\ \vdots \\ \left(q^{i,j}_l\right)^*\left(q^{j,i}_l\right)\end{bmatrix}  = \begin{bmatrix}\theta_{{i,q}_1} \\ \vdots \\ \theta_{{i,q}_l} \end{bmatrix} \in \mathbb{R}^{4l},\label{eq:quat1}
    \end{align}
where the product shown here can be understood as a quaternion multiplication which acts on each quaternion in the vector independently.\\

The last step needed for the soft consensus is to correctly define a tuning constant $\mu$, which acts as a feedback gain on the quaternion consensus term in Eq.~\eqref{eq:quat1}.
To preserve the normalization of the quaternion, we define the quaternion scaling function,{
\begin{equation}\label{eq:quat_scale}
    q_s(q, \mu) = \begin{bmatrix} \sqrt{1 - \mu^2\sum_{k=1}^3 q_{k}^2}\\
    \mu q_{1}\\
    \mu q_{2}\\
    \mu q_{3}
    \end{bmatrix},
\end{equation}
}
and with this function the quaternion error is scaled appropriately as,
\begin{equation}\label{eq:distributed_soft_cons_attitude}
    \phi_{\text{full}, i, q}^\mu := \begin{bmatrix} q_s(\theta_{{i,q}_1}, \mu_{i, q}) \\ \vdots \\ q_s(\theta_{{i,q}_l}, \mu_{i, q})\end{bmatrix} \in \mathbb{R}^{4l}.
\end{equation}
As a consistency check, as $\mu\to0$, the attitude approaches $[1~0~0~0]$, which denotes an error of zero.
Finally, as for the bias and position vectors, the last step is to reduce this quantity to the appropriate size by only keeping the entries which are in the set ${V}_i$, yielding $\phi_{i, q} \in \mathbb{R}^{4(k_i+1)}$.

\textit{Combining the Soft Consensus Terms}\\
Now that the three soft consensus terms are computed, they need to be applied to the estimate of satellite $i$. 
To apply the correction term on the attitude, the current estimate of the attitude is multiplied (in the quaternion sense) with the attitude component of the soft consensus term $\phi_{i, q}$ as, 
\begin{align}
    \overline{\hat q_i^{++}} &= \overline{\hat q_i^+}\phi_{i, q} = \begin{bmatrix}\hat q^+_{i, j_{i,1}}\phi_{{i, q}_1} \\ \vdots \\ \hat q^+_{i, j_{i,k_i+1}}\phi_{{i, q}_{k_i+1}}\end{bmatrix}.
    \end{align}
The bias estimate is corrected by adding the bias soft consensus term as,
\begin{align}
      \overrightarrow{\hat b_i^{++}} &= \overrightarrow{\hat b_i^+} + \phi_{i, b}.
\end{align}
Finally, the position soft consensus term is added to Eq.~\eqref{eq:distributed_kf_position_def}
\begin{align}
    \overline{\hat r_i^{++}}& = \overline{\hat r_i^+} + \phi_{i, r},
\end{align}
 and then transformed to the dual quaternion position, as
\begin{align}
    \overrightarrow{\hat q^{++}_i} &= 
       \begin{bmatrix} 
         \hat q^{++}_{i, j_{i,1}} \\ 
         \frac{1}{2} 
         \begin{bmatrix} 
           0 \\ \hat r^{++}_{i, j_{i,1}} 
         \end{bmatrix} \hat q^{++}_{i, j_{i,1}}  \\ 
         \vdots \\  
         \hat q^{++}_{i, j_{i,k_i+1}} \\ 
         \frac{1}{2} 
         \begin{bmatrix} 
           0 \\ \hat r^{++}_{i, j_{i,k_i+1}} 
         \end{bmatrix} 
         \hat q^{++}_{i, j_{i,k_i+1}}
       \end{bmatrix}.
\end{align}

\subsubsection{Description of the Hard Consensus Term}
\label{sec:descr-hard-cons}
The purpose of the hard consensus term is to fuse measurement data and error covariance matrices between neighbouring agents.
We adapt the method proposed in~\cite{1583238} to the dual quaternion setting. 
The main challenge is computing the distributed versions of the Kalman information filter quantities in~\eqref{eq:y_S_inf_form_def}, as each satellite has a different number of estimated states and so the measurement matrices have different sizes for each node and also incorporate different measurements. 
 As such, each satellite needs to prepare a measurement matrix and vector compatible with each of its neighbours.

Consider the measurement matrix of satellite $i$ as defined in Eq.~\eqref{eq:distributed_H}. 
Then, $H_{i, i} = H_{i} \in \mathbb{R}^{8(k_i+1)\times 14(k_i+1)}$.
We will define the matrix $H_{i, k}$ as just a reconfiguration of $H_{i,i}$ to match the different state of satellite $k$, as shown in Appendix~\ref{sec:full-jacobian-hard} Eqs.~\eqref{eq:H_hard_1}-\eqref{eq:H_hard_2}.
Because some states are not estimated by satellite $i$ but are estimated by its neighbour satellite $k$, it may be the case that depending on the neighbourhood, some of the measurement Jacobians satisfy $H_{i, n_{k,1}}^{n_{k,1}} = \mathbf{0}_{8\times 8}$.

A similar procedure can now be done for the measurement itself.
As for the measurement Jacobian, each satellite needs to prepare the measurement and convert it to a vector compatible with the state estimate of the neighbouring satellites. 
The measurement of satellite $i$ is defined as in Eq.~\eqref{eq:measurement_vec_distr}. 
Then, the vector $v_o \in \mathbb{R}^6$ (being the $o$-th entry of $z_{i, k}$) is equal to the vector $c_l \in \mathbb{R}^6$ (being the $l$-th entry of $z_i$) if $j_{i,l} = j_{k, o}$. Otherwise the entry is $\mathbf{0}_{6\times 1}$.
Explicitly, for $k \in {V}_i$,
\begin{align}
    z_{i,k} = \begin{bmatrix}v_1 \\ \vdots \\ v_{k_k+1}\end{bmatrix} \quad \text{where }\begin{cases}
    v_o = c_l, \qquad &\text{if } j_{i,l} = j_{k, o}\\
    v_o = \mathbf{0}_{6\times 1}, \qquad &\text{otherwise}
    \end{cases}.
\end{align}
In the running example described in Fig~\ref{fig:swarm-network}, the sensitivity matrices and measurement vectors of satellite 4 for the estimates of itself and its neighbours 1 and 3 are shown in Eqns~\eqref{eq:egH1}, \eqref{eq:egH2} and \eqref{eq:egH3}.
\begin{subequations}
 \label{eq:egH}
  \begin{align}
    \overline{H_{4, 1}} &= \begin{bmatrix}
            \overline{H^1_{4, 1}} &  \mathbf{0}_{6 \times 6} & \mathbf{0}_{6 \times 6} & \mathbf{0}_{6 \times 6} & \overline{H^4_{4, 1}} & \mathbf{0}_{6 \times 6}\\
            \mathbf{0}_{6 \times 6} & \mathbf{0}_{6\times 6} & \mathbf{0}_{6\times 6} & \mathbf{0}_{6 \times 6} & \mathbf{0}_{6\times 6} & \mathbf{0}_{6 \times 6}\\
            \mathbf{0}_{6 \times 6} & \mathbf{0}_{6\times 6} & \mathbf{0}_{6\times 6} & \mathbf{0}_{6 \times 6} & \mathbf{I}_{6\times 6} & \mathbf{0}_{6 \times 6}
            \end{bmatrix}\in \mathbb{R}^{18\times 36} 
            &
            z_{4,1} = \begin{bmatrix}\overline{q_{B_4/B_1, m}} \\ \overline{r_{B_1/B_4, m}^{B_4}} \\ \mathbf{0}_{6\times 1}\\ \overline{\dq{\hat q}_{4,4/I}\dq{q}_{B_4/I, m}}\end{bmatrix}\label{eq:egH1}\\
    \overline{H_{4, 3}} &= \begin{bmatrix}
            \mathbf{0}_{6 \times 6} & \mathbf{0}_{6\times 6} & \mathbf{0}_{6\times 6} & \mathbf{0}_{6 \times 6} & \mathbf{0}_{6 \times 6} & \mathbf{0}_{6\times 6} & \mathbf{0}_{6\times 6} & \mathbf{0}_{6 \times 6}\\
            \mathbf{0}_{6 \times 6} & \mathbf{0}_{6\times 6} & \overline{H^3_{4, 3}} &  \mathbf{0}_{6 \times 6} & \overline{H^4_{4, 3}} & \mathbf{0}_{6 \times 6} & \mathbf{0}_{6 \times 6} & \mathbf{0}_{6\times 6}\\
            \mathbf{0}_{6 \times 6} & \mathbf{0}_{6\times 6} & \mathbf{0}_{6 \times 6} & \mathbf{0}_{6\times 6} & \mathbf{I}_{6\times 6} & \mathbf{0}_{6 \times 6} & \mathbf{0}_{6 \times 6} & \mathbf{0}_{6\times 6}\\
            \mathbf{0}_{6 \times 6} & \mathbf{0}_{6\times 6} & \mathbf{0}_{6\times 6} & \mathbf{0}_{6 \times 6} & \mathbf{0}_{6 \times 6} & \mathbf{0}_{6\times 6} & \mathbf{0}_{6\times 6} & \mathbf{0}_{6 \times 6}
            \end{bmatrix}\in \mathbb{R}^{24\times 48}
            &
            z_{4,3} = \begin{bmatrix}\mathbf{0}_{6\times 1} \\ \overline{q_{B_4/B_3, m}} \\ \overline{r_{B_3/B_4, m}^{B_4}} \\ \overline{\dq{\hat q}_{4,4/I}\dq{q}_{B_4/I, m}}\\ \mathbf{0}_{6\times 1}\end{bmatrix}\label{eq:egH2}\\
    \overline{H_{4, 4}} &= \begin{bmatrix}
            \overline{H^1_{4, 1}} &  \mathbf{0}_{6 \times 6} &  \mathbf{0}_{6 \times 6} &  \mathbf{0}_{6 \times 6} & \overline{H^1_{4, 1}} & \mathbf{0}_{6 \times 6}\\
            \mathbf{0}_{6 \times 6} &  \mathbf{0}_{6 \times 6} &\overline{H^3_{4, 3}} &  \mathbf{0}_{6 \times 6} & \overline{H^4_{4, 3}} & \mathbf{0}_{6 \times 6}\\
            \mathbf{0}_{6 \times 6} & \mathbf{0}_{6\times 6} & \mathbf{0}_{6 \times 6} & \mathbf{0}_{6\times 6} & \mathbf{I}_{6\times 6} & \mathbf{0}_{6 \times 6}
            \end{bmatrix}\in \mathbb{R}^{18\times 36}
            &
            z_{4,4} = \begin{bmatrix}\overline{q_{B_4/B_1, m}} \\ \overline{r_{B_1/B_4, m}^{B_4}} \\ \overline{q_{B_4/B_3, m}} \\ \overline{r_{B_3/B_4, m}^{B_4}} \\ \overline{\dq{\hat q}_{4,4/I}\dq{q}_{B_4/I, m}}\\ \mathbf{0}_{6\times 1}\end{bmatrix}\label{eq:egH3}
  \end{align}
\end{subequations}

Now that the measurement matrices and vectors are compatible between neighbours, we can formalize the consensus step. 
The only equations which change in comparison to the distributed Kalman filter without a consensus step are Eqs.~\eqref{eq:meas_cov_distr1}-\eqref{eq:meas_cov_distr3}. 
Each $i\in V$ satellite prepares the measurement covariance matrix and the measurement for its neighbour $k\in N_i$,
\begin{equation}
    \begin{aligned}
    \overline{U_{i, k}} &= \overline{H_{i, k}^T }R_{i, k}^{-1}\overline{H_{i, k}} \\ 
    U_{i, k} &= H_{i, k}^T  R_{i, k}^{-1}H_{i, k}\\
    u_{i, k} &= \overline{H_{i,k}^T } R_{i, k}^{-1}\overline{z_{i, k}},
    \end{aligned}
\end{equation}
and sends a message to its neighbour $k$,
\begin{equation}
    \text{msg}_{i,k} = \left(\overline{U_{i,k}}, U_{i,k}, u_{i,k}\right).
\end{equation}
Next, each satellite $i\in V$ calculates the aggregated covariance matrix $S_i$ and the aggregated sensor data $y_i$ as,
\begin{equation}\label{eq:hard_cons_aggregation}
    y_i = \sum_{j\in {V}_i} u_{j, i}\qquad
    S_i = \sum_{j\in {V}_i} U_{j, i}\qquad
    \overline{S_i} = \sum_{j\in {V}_i} \overline{U_{j, i}}.
\end{equation}
Each satellite $i\in V$ then computes the state update~\eqref{eq:M_matrix_distributed} as, 
\begin{equation}
    \Delta x_i = 
    \begin{bmatrix}
    \Delta \overline{\dq{\delta \hat{q}}_{j_{i,1}/I}^T} &
    \Delta \overline{\dq{\delta \hat{b}}_{\omega_{j_{i,1}}}^T} &
    \cdots&
    \Delta \overline{\dq{\delta \hat{q}}_{j_{i,k_i+1}/I}^T} &
    \Delta \overline{\dq{\delta \hat{b}}_{\omega_{j_{i,k_i+1}}}^T} 
    \end{bmatrix}^T
    = M_i\left(\overline{y_i} - \overline{\left(S_i \mathbf{1}\right)}\right),
\end{equation}
where
\begin{equation}\label{eq:M_matrix_distributed_2}
    M_i = \left(P_i^{-1} + \overline{S_i}\right)^{-1}.
\end{equation}
Finally, each agent $i\in V$ updates the estimates of the pose and bias using Eqs.~\eqref{eq:extend_to_R8_dist} and \eqref{eq:dist_meas_update_estimated}.
\subsection{Leader-Follower Consensus}
\label{sec:lead-foll-cons}
In a satellite fleet where absolute pose measurements are difficult to obtain or would require an expensive sensor (i.e., a star-tracker), it might be interesting to only provide pose measurements to a subset of the satellite fleet. 
The remaining satellites would then only have relative pose measurements with their neighbours. 

In this configuration, we define a set of \emph{leaders} $L \subset {V}$ which includes all satellites that have access to absolute pose measurements. 
We then denote $F = V\setminus L$ as the \emph{followers}.

The removal of the pose measurement for the followers only changes their respective measurement update of the Kalman filter. 
A straightforward modification consists of removing the absolute pose measurement from the measurement vector and remove the corresponding row of the measurement matrix in~\eqref{eq:distributed_H} for all satellites $i \in F$.
Empirically, this solution does not yield good results if the set $L$ is small. 
Hence, a different approach was chosen for the estimation of the followers.
\subsubsection{Calculation of Absolute Pose Measurements}
To address this issue, we propose a method whereby an absolute pose measurement is constructed from relative pose measurements and absolute pose estimates. 
At each timestep, satellite $i\in V$ estimates its own pose and the pose of its neighbours $k\in N_i$. 
Additionally, it measures the relative pose of the neighbouring satellites in its own body frame. 
Satellite $i$ can thus compute $|N_i| = k_i$ absolute pose measurements for $k=1,\dots,k_i$ through
\begin{align}
    q_{B_i/I, m_k} &= \hat q^-_{i, n_{i, k}}q_{B_{n_i, k}/B_i, m}^*\label{eq:abs_pose_calc_1}\\
    r_{B_i/I, m_k}^I &= 2 \hat q^-_{i, n_{i, k}} \hat p^-_{i, n_{i, k}} - q_{B_{n_i, k}/B_i, m} r_{B_{n_i, k}/B_i, m}^{B_i} q_{B_{n_i, k}/B_i, m}^*.\label{eq:abs_pose_calc_2}
    \end{align}
Informally, to get an absolute attitude measurement of satellite $i$, we remove the relative attitude measurement from satellite $i$ to satellite $k$ from the estimated attitude of satellite $k$ computed by satellite $i$. The same is done for the position measurement.
To fuse the various position estimates, we can average them as,
\begin{equation}
    \overline{r_{B_i/I, m}^I} = \frac{1}{k_i}\sum\limits_{k \in N_i} \overline{r_{B_i/I, m_k}^I}.
\end{equation}
The notion of the `average' of the attitude quaternion is proposed by~\cite{quaternion_averaging} and is computed as the minimizer of the linear program,
\begin{align}
    \overline q_{B_i/I, m} &= \left\{
                             \begin{array}{ll}
                               \mathrm{argmax} & \mathrm{Tr}\left(A(q_{B_i/I, m})C^T \right)\\
                               \mathrm{subject~to} & C  = \sum\limits_{k\in N_i} A(q_{B_i/I, m_k}),
                             \end{array}\right.
\end{align}
where $A(q)$ is the rotation matrix of quaternion $q$.
The dual quaternion for the pose measurement is then constructed from $\overline{r_{B_i/I, m}^I}$ and $\overline q_{B_i/I, m}$, and no further changes to the algorithm are necessary.

\subsubsection{Stubborn Leaders}
\label{sec:dist_variant}
In the leader-follower setting, it is quite clear that leaders always have more information than followers due to their absolute pose measurement. 
In many cases, a follower may only be able to improve the estimate of a leader if it has immediate knowledge of another leader satellite, i.e., if the follower is connected to two or more leaders. 
Notwithstanding pathological cases, this occurs in a graph when a large percentage of satellites are leaders. 
Conversely, as the percentage of leaders becomes small, the follower satellites may very well reduce the estimator performance of the leader, which in turn may deteriorate the performance of the whole fleet. 
In this section, we describe the notion of a `stubborn' leader which does not perform consensus with followers, which, as we will see in \S\ref{sec:lead-foll-scen}, provides better average performance when the percentage of leaders is low.
Only two changes need to be made.
In the case of soft consensus, it sufffices to set
$    \mu_{i, r} = \mu_{i, q} = \mu_{i, b} = 0 $ if $i\in L$
for leaders to ignore the follower estimates. 
In the case of hard consensus, Eq.~\eqref{eq:hard_cons_aggregation} must be replaced with,
\begin{align}
    y_i &= \begin{cases} u_{i,i}, & \text{if } i \in L\\
                        \sum\limits_{j\in {V}_i}u_{j, i}, & \text{otherwise,}\end{cases}~
    S_i = \begin{cases} U_{i,i}, & \text{if } i \in L\\
                        \sum\limits_{j\in {V}_i}U_{j, i}, & \text{otherwise,}\end{cases}~
    \overline{S_i} = \begin{cases} \overline{U_{i,i}}, & \text{if } i \in L\\
                        \sum\limits_{j\in {V}_i}\overline{U_{j,i}}, & \text{otherwise.}\end{cases}
    \end{align}
With these modifications, the leaders will not be influenced by the follower estimates.

\section{Numerical Examples}
\label{sec:numerical-examples}

In this section, we present three numerical experiments that demonstrate the advantages of the distributed DQ-MEKF.
First, in \S\ref{sec:comp-single-satell} we compare the DQ-MEKF and its distributed counterpart.
Next, in \S\ref{sec:coop-swarm-around} we showcase the filter in a scenario where a fleet of satellites swarm around an asteroid.
Finally, in \S\ref{sec:lead-foll-scen}, we examine the performance of the distributed DQ-MEKF in the leader-follower scenario where only a subset of the satellites have an absolute pose measurement.
We provide an open-source repository~\cite{binzcode} containing the code to reproduce the numerical examples in this section.

\subsection{Comparison of Cooperative vs Non-Cooperative Fleet}
\label{sec:comp-single-satell}

\begin{table}[ht]
\centering

\begin{tabular}[t]{lc}
\hline
Matrix & Value\\
\hline
$\overline{Q}_{b_{\omega}}$ & $\frac{10^{-3}}{\text{SNR}^2}\mathbf{I}_{3\times 3}$\\
$\overline{Q}_{b_v}$ & $\frac{10^{-1}}{\text{SNR}^2} \mathbf{I}_{3\times 3}$\\
$\overline{Q}_{\omega},\overline{Q}_v$ & $\mathbf{0}_{3\times 3}$\\
$R_{6\times 6}$ & $\mathrm{blkdiag}\{\mathrm{std}_q^2 \mathbf{I}_{3\times3} ,\mathrm{std}_r^2\mathbf{I}_{3\times 3}\}$\\
$P_{12\times 12}(0)$ & $\mathrm{blkdiag}\{10^{-1}\mathbf{I}_{6\times 6},  10^{-2}\mathbf{I}_{6\times 6}\}$\\
$P_{15\times 15}(0)$ & $\mathrm{blkdiag}\{10^{-1}\mathbf{I}_{6\times 6},  10^{-2}\mathbf{I}_{9\times 9}\}$\\
\hline
\end{tabular}
\caption{Simulation parameters for the comparison between the DQ-MEKF and DDQ-MEKF.}
\label{tab:covariances_single_sat}
\end{table}

We aim to compare the performance of Algorithm~\ref{alg:ddq-mekf} on a fleet of satellites with that same fleet of satellites running the DQ-MEKF algorithm from~\cite{filipe2015extendedjgcd}, as this essentially quantifies the improvement of performance when distributing the algorithm.

In particular, we compare the performance of the DQ-MEKF, Algorithm~\ref{alg:ddq-mekf} with both hard and soft consensus, and Algorithm~\ref{alg:ddq-mekf} with soft consensus only.
We assume absolute pose measurements for all satellites, and in each case, the respective filter was initialized with the parameters in Table~\ref{tab:covariances_single_sat}.
In the distributed case, for each satellite $i\in V$, the covariance matrices for the various noises are direct-summed with a copy for each neighbour $j\in V_i$.
Whenever soft consensus is employed, we set $ \mu_{i, r} = \mu_{i, q} = \mu_{i, b} = 1/|J_i|$ where $|J_i|$ is the number of neighbours of $i$, including itself.

Each simulation is initialized with a given signal-to-noise ratio (SNR), and we do a simulation sweep over values of the SNR from 5 to $10^5$.
In Table~\ref{tab:covariances_single_sat},  $\mathrm{std}_q$ and $\mathrm{std}_r$ are the standard deviations of the Gaussian noise applied to the attitude and position respectively, which are calculated from the given SNR.
For each algorithm and each SNR, 80 such simulations were run at 20Hz for 60s.
Finally, the graphs defining the connectivity of the fleet of satellites were generated randomly, with the probability of an edge existing between $i,j\in V$ being set to 0.5. 
Disconnected graphs generated in this manner were discarded and re-initialized.

\begin{figure}[!htb]
  \centering
    \includegraphics[width=\textwidth]{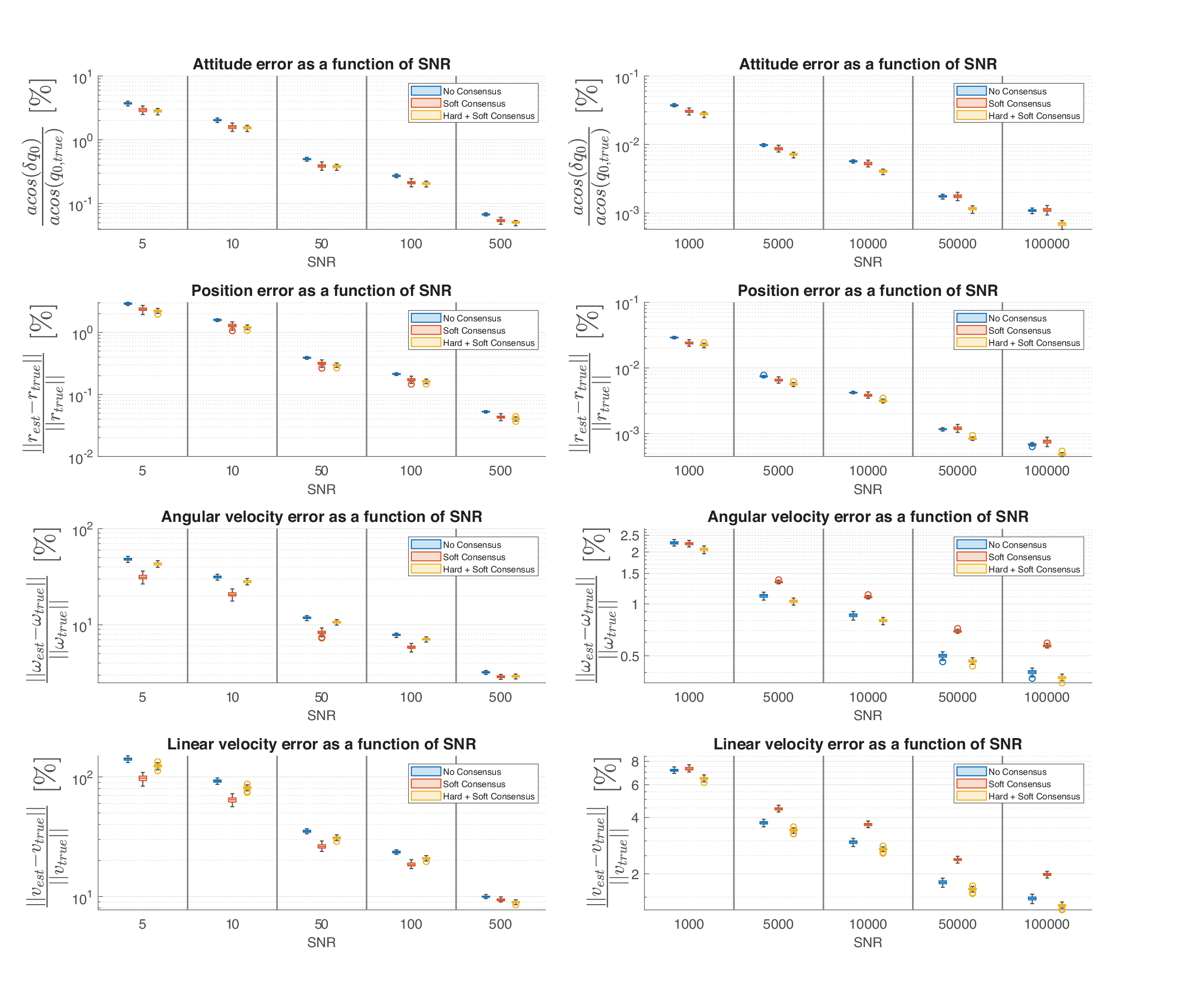}
    \caption{Comparison of the different consensus filters for a fleet of ten satellites. Blue: the baseline non-cooperative DQ-MEKF. Red: the distributed DQ-MEKF is shown with the soft consensus term. Yellow: the distributed DQ-MEKF with soft and hard consensus.}
    \label{fig:comparison_cooperating_changing_graph}
\end{figure}

The results of the parameter sweep are depicted in Fig.~\ref{fig:comparison_cooperating_changing_graph}, where we plot the interquartile ranges of attitude, position, angular velocity and linear velocity estimator error over the simulations for varying SNRs.
{For small SNRs, soft consensus is is superior to the non-cooperative case and performs similarly to the case with both soft and hard consensus; however, for higher SNRs, soft consensus provides an inferior velocity estimate compared to the other filters.}
However, the DDQ-MEKF with both hard and soft consensus outperforms the non-cooperative case in every metric over all SNRs.

\subsection{Distributed Swarming Around an Asteroid}
\label{sec:coop-swarm-around}

\begin{table}[ht]
\centering
\begin{tabular}[t]{l c c c c}
\hline
Parameter $(\forall i\in V)$ &Value\\
\hline
    Attitude Noise Cov & $ 2.79\times 10^{-7} \mathbf{I}_{4\times 4}$ \\
    Position Noise Cov & $8.55\times 10^{-4} \mathbf{I}_{3\times 3} $\\
    $\overline Q_{\dq{\omega}_i}$ & blkdiag$\{0\times \mathbf{I}_{6\times6}\}$\\
    $\overline Q_{\dq{b_\omega}_i}$ & blkdiag$\{10^{-6}\mathbf{I}_{3\times 3},  10^{-4}\mathbf{I}_{3\times 3}\}$\\
    $P_i(0)$ & blkdiag$\{10^{-1}\mathbf{I}_{6\times 6}, 10^{-2}\mathbf{I}_{6\times 6}\}$ \\
    $\mu_{i, q}, \mu_{i, r},\mu_{i, b}$ & $1/|J_i|$\\
    $Q_{LQR, i}$ & $1\times 10^{-1} \mathbf{I}_{12\times 12}$\\
    $R_{LQR, i}$ & $1\times 10^{-1} \mathbf{I}_{6\times 6}$\\
\hline
\end{tabular}
\caption{Settings of the distributed DQ-MEKF for the distributed swarming example.}
\label{tab:settings_case_study}
\end{table}%

\begin{figure}
  \centering
  \includegraphics[width=0.5\columnwidth]{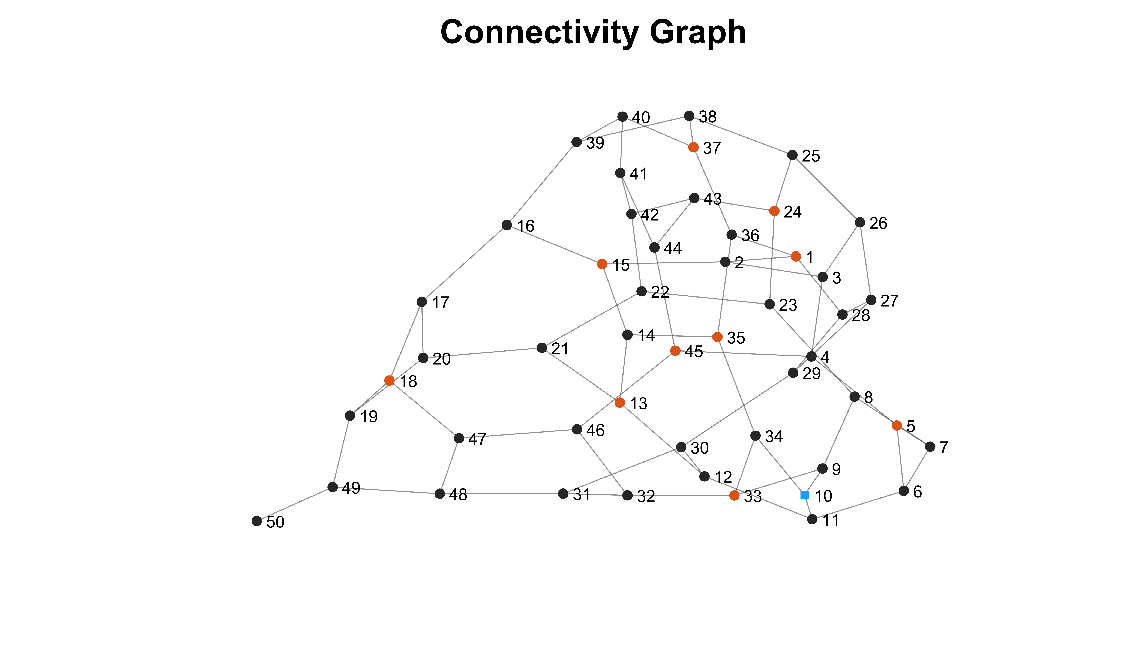}
  \caption{Connectivity graph for the asteroid coverage scenario. Red nodes depict leaders, and the blue square node depicts satellite 10.}
  \label{fig:swarm-network}
\end{figure}

In this example, we consider the performance of Algorithm~\ref{alg:ddq-mekf} (with both hard and soft consensus) in a scenario where 50 satellites arrange themselves around an asteroid, while maintaining an attitude representative of pointing a scientific payload at the asteroid.
The initial conditions of the satellites are arranged on a plane 40m from the asteroid, and a Fibonacci lattice~\cite{fibonacci_lattice} is used to evenly space the target positions of the satellites on a sphere of radius 25m centered around the asteroid.
The reference attitude is constructed from the current estimate of the satellite and the position of the asteroid, and an LQR controller is used to maneuver the fleet to their final positions and attitudes over \textasciitilde 200s.
During the maneuver, the satellites perform Algorithm~\ref{alg:ddq-mekf} across the connectivity graph seen in Fig.~\ref{fig:swarm-network}, with a sampling rate of 20 Hz, and initialized with the parameters in Table~\ref{tab:settings_case_study}.

\begin{figure}
  \centering
  \includegraphics[width=0.5\columnwidth]{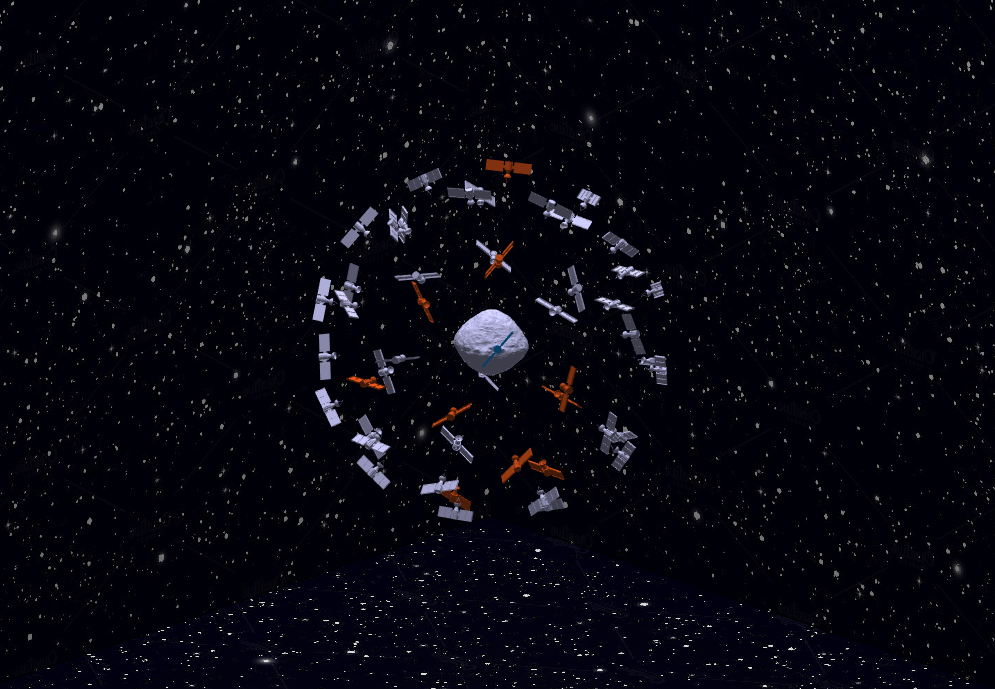}
  \caption{Target formation of the cooperative asteroid swarming scenario. Satellites depicted in red are the leaders.}
  \label{fig:asteroid_pos}
\end{figure}

\begin{figure*}
  \centering
  \includegraphics[width=0.35\textwidth]{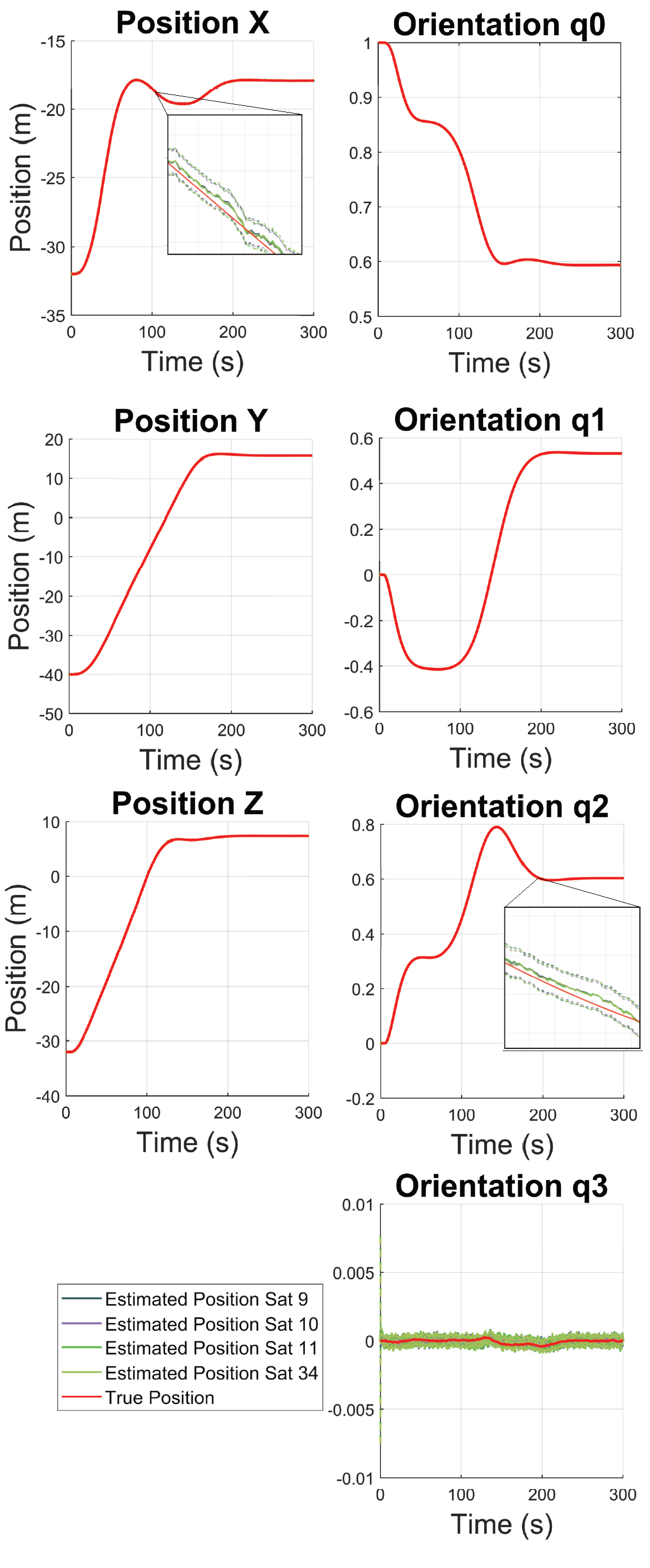}
  \caption{Attitude and position of satellite 10 in the asteroid coverage scenario.}
  \label{fig:asteroid-states}
\end{figure*}

\begin{figure*}
  \centering
  \includegraphics[width=0.7\textwidth]{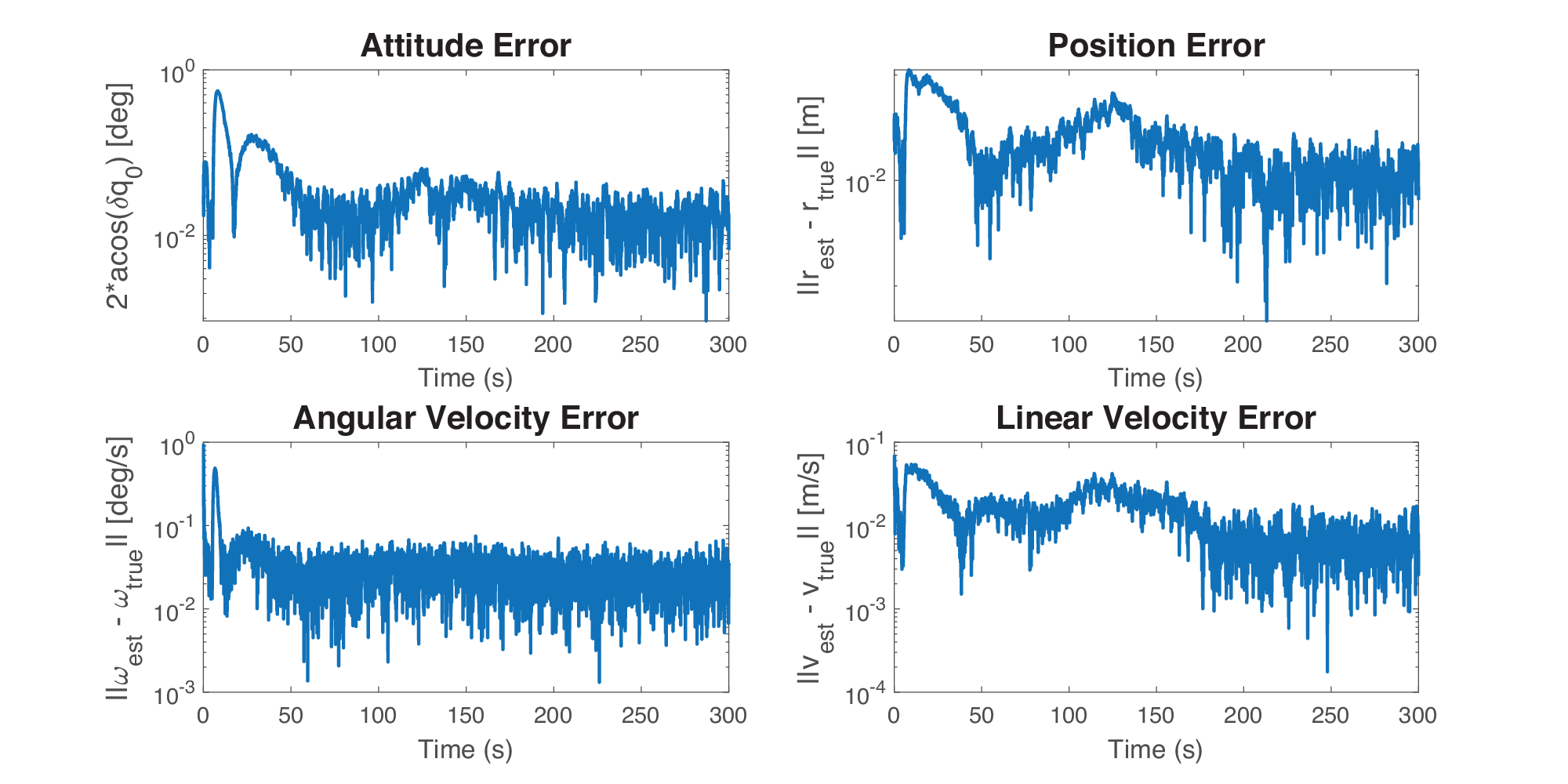}
  \caption{RMS error of the attitude and position and their derivatives of satellite 10 in the asteroid coverage scenario.}
  \label{fig:asteroid-rms}
\end{figure*}

A representative satellite with node label 10 {is selected at random}, and is coloured in blue in Figs.~\ref{fig:asteroid_pos} and~\ref{fig:swarm-network}.
{Leader nodes are depicted in red.}
The position and attitude (in quaternion form) of satellite 10 throughout the maneuver is shown in Fig.~\ref{fig:asteroid-states}
, along with the estimates of satellite 10 as computed by Algorithm~\ref{alg:ddq-mekf} by satellite 10 itself, and its neighbours.
The RMS errors of the attitude, position, angular velocity, and linear velocity are shown in Fig.~\ref{fig:asteroid-rms}.
As one can see, the error decreases over time until it reaches a steady-state, showing that the estimator remains stable during the maneuver.
A supplementary video showing the swarming behaviour of entire the satellite fleet, as well as real time pose data from satellite 10, may be found at~\cite{ourvideo}.

\subsection{Leader-Follower Scenario}
\label{sec:lead-foll-scen}

\begin{figure}[!htb]
  \centering
    \includegraphics[width=0.5\textwidth]{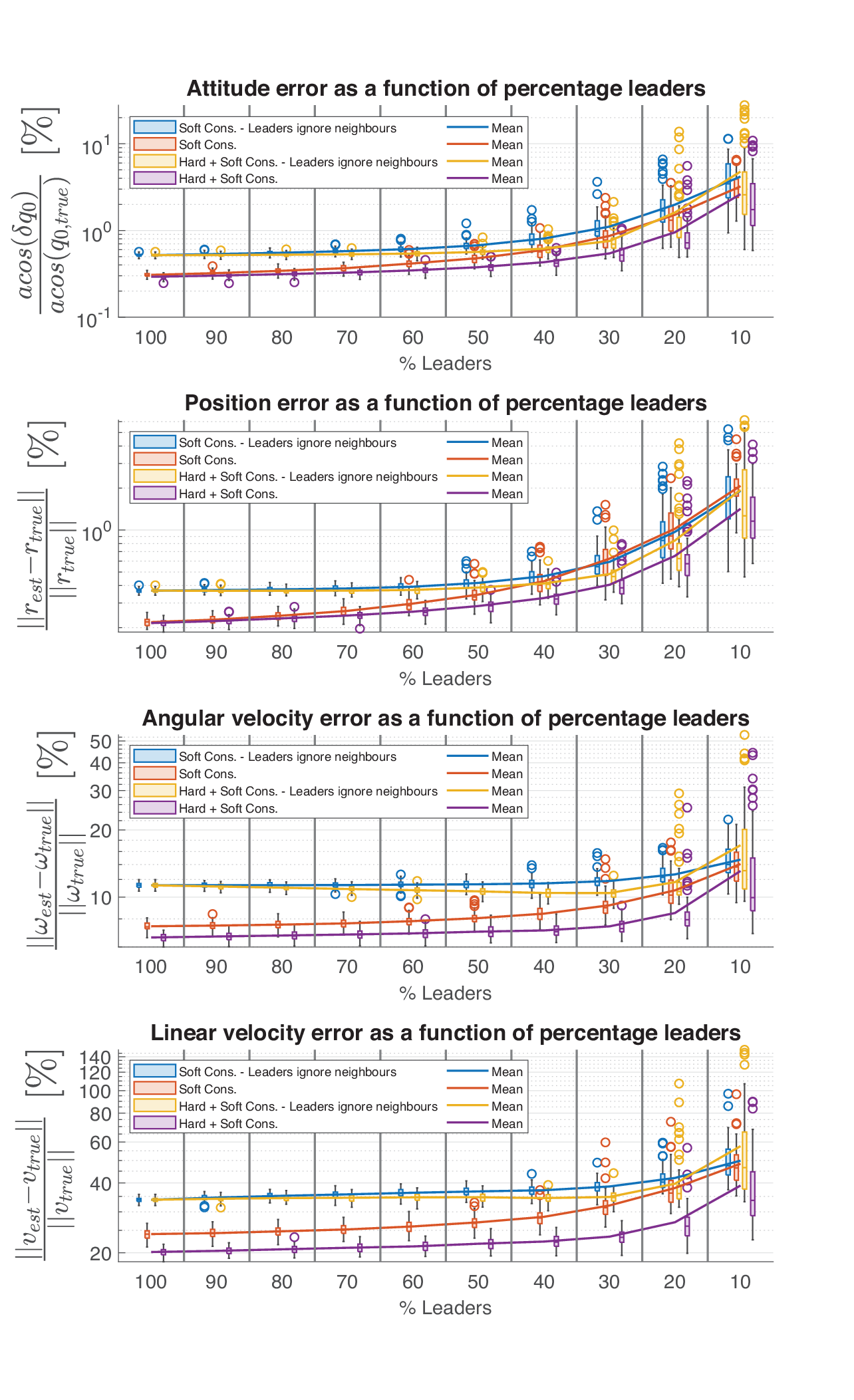}
    \caption{Comparison of the four variations of the  DDQ-MEKF by varying the \% of leaders over a random graph. The fleet has ten satellites, and has a fixed SNR of 1000.}
    \label{fig:comparison_leaders_changing_graph}
\end{figure}

In this section, we provide a set of numerical experiments to investigate the performance of Algorithm~\ref{alg:ddq-mekf} with leader-follower consensus as described in \S\ref{sec:lead-foll-cons}.
Explicitly, we consider four cases:
with hard and soft consensus, with soft consensus only,  with hard and soft consensus and subborn leaders, and with soft consensus and stubborn leaders.
In each case, only the leaders have access to absolute pose measurements, and followers create synthetic pose measurements using Eqn.~\eqref{eq:abs_pose_calc_1}.

For each experiment, the SNR was fixed to 1000, and the filters were initialized with the parameters in Table~\ref{tab:covariances_single_sat}.
We varied the percentage of leaders from 10\%-100\%, and leaders were selected {uniformly} at random in the graph {in order to understand how the performance varies with different sets of leaders}.
The graph itself was randomized as in \S\ref{sec:comp-single-satell}.
 {For each instance of leader percentage,  the simulation was run 80 times (randomizing the graph and leaders each time), and each simulation was run for a duration of 60s at 20Hz (1200 iterations).}
{In order to understand the steady-state error of each simulation,}
 the error values in Fig.~\ref{fig:comparison_leaders_changing_graph} are computed using the last 600s of each simulation to discard transient behaviour.

We can see in Fig.~\ref{fig:comparison_leaders_changing_graph} that algorithms with stubborn leaders perform worse than with non-stubborn leaders, and overall the algorithms perform better with larger percentages of leaders.
This is expected.
When leaders ignore their neighbours and all satellites are leaders, this is equivalent to the case in \S\ref{sec:comp-single-satell} with non-cooperating satellites, where the cooperating fleet has superior performance to that of the non-cooperating fleet. 
The soft consensus case also always does worse than the combination of hard and soft consensus for a SNR of 1000, which is in agreement with \S\ref{sec:coop-swarm-around}.

As expected, as the percentage of leaders decreases, the performance of the entire fleet also decreases.
It can also be seen that the case of stubborn leaders doesn't suffer as much performance degradation as the case where all satellites perform consensus. 
This is due to the fact that the followers cannot have a better estimate of the leader if the number of leaders is very small.
A final observation is that the performance of the fleet doesn't degrade until approximately less than 50\% of fleet are leaders.
Since the graphs and leaders are randomly generated, this is the point at which a follower has a high likelihood of only being connected to one leader, or to no leaders at all. 
Therefore, its performance is limited by the quality of measurements it has access to, which in turn is determined in part by the quality of the estimates of its neighbours.
If this neighbourhood doesn't contain leaders, then the quality of the corresponding satellite's estimate reflects the fact that the absolute pose measurements provided by the leaders is degraded as it propagates across the network.

\section{Conclusion}
\label{sec:conclusion}

In this paper, we have proposed, derived and validated a distributed version of the dual-quaternion multiplicative extended Kalman filter as studied by~\cite{filipe2015extendedjgcd}.
The relevant algebraic quantities for distributing the filter were derived, and the concepts of hard and soft consensus were introduced to improve measurement and estimate fusion between neighbouring satellites in the network.
We also described the filter in a leader-follower scenario, where only leaders in the fleet have access to an absolute pose measurement.

Our numerical experiments show that distributing the DQ-MEKF over a network of satellites provides a substantial performance increase across a variety of signal-to-noise ratios over the non-cooperative scenario when both hard and soft consensus is employed alongside Algorithm~\ref{alg:ddq-mekf}.
Furthermore, we established numerically that the filter is sufficiently stable and provides an estimate of sufficient fidelity that it can be succesfully used in a dynamic scenario where satellites are cooperating to swarm around an asteroid.
This experiment also demonstated the scalability of the algorithm over a larger fleet.
Our final numerical experiment showed that we can remove expensive absolute pose measurement sensors, such as star-trackers, from a large part of the fleet with no significant performance degradation.
{We leave the question of optimally choosing which satellites to select as leaders (given a network) for future research.}

{Additional} future work may include validation of the proposed methods on a hardware testbed involving multiple satellite simulators, or adapting the filter to pose control of other vehicles such as quadrotors.
{Considering correlations between neighbouring pose estimates in the consensus algorithms may yield further improvements in performance}.
One remaining question not fully addressed in this paper is a detailed analysis  of the effect of the network topology on the performance of the algorithm -- our chosen network topologies were randomly generated, and as such the average effect of the networks were studied, but the question of optimizing specific networks for the filter remains.
{Aside from the dynamical analysis contained in this paper, this line of work will involve topics in combinatorial optimization.}

\appendix
\section{Appendix}
\label{sec:appendix-2}

\subsection{Jacobians of Measurement Models}
\label{sec:jacob-meas-models}
Here, we derive the Jacobians of the measurement functions in~\eqref{eq:measurement_vec_distr} for various measurement models.
For the sake of reducing notational burden, the time dependence of the estimates is neglected in the derivations below. 
However, we ask the reader to keep in mind that the estimates and measurements are also time-dependent.

\subsubsection{Jacobian of the Absolute Pose Measurement}
As can already be seen in Equation \eqref{eq:measurement_vec_distr}, the measurement is multiplied with the conjugate of the estimate, directly yielding the error dual quaternion. 
Hence, the measurement matrix $H$ is simply,
\begin{equation}\label{eq:distr_abs_meas_H_matrix}
H_{abs} = 
\begin{bmatrix}
  \mathbf{I}_{6\times 6}& \mathbf{0}_{6\times 6}
\end{bmatrix}.
\end{equation}
\\

\subsubsection{Jacobian of the Relative Attitude Measurement}
The relative attitude $q_r$ is the real part of $\dq{q}_{B_k/B_i} = q_{r, B_k/B_i} + \epsilon q_{d, B_k/B_I}$. 
To reduce notational burden, we will denote $q_{r, B_k/B_i} := q_{k/i}$ and $q_{d, B_k/B_I} := p_{k/i}$. 
In this way, $q_{k/i}$ is the quaternion describing the rotation from satellite $i$ to satellite $k$, and $p_{k/i}$ denotes the position quaternion from satellite $i$ to $k$. 

We define the error quaternion as $\dq{\delta q}_{B/I} = \dq{\hat q}_{B/I}^* \dq{q}_{B/I}$, where $B$ is the body frame of the satellite. 
Through the definition of the dual quaternion, it immediately follows that $\delta q_{B/I} = \hat q_{B/I}^* q_{B/I}$. 
 Dropping $B$ from the notation as above yields,
\begin{align}
    \delta q_{i/I} &= \hat q_{i/I}^* q_{i/I} \\
    \delta q_{k/I} &= \hat q_{k/I}^* q_{k/I},
\end{align}
 which allows us to compute $h_{q_{k/i}}(\cdot)$ as,
 \begin{align}
   q_{k/i} &= q_{i/I}^*q_{k/I}\\
           &= (\hat q_{i/I}\delta q_{i/I})^*(\hat q_{k/I}\delta q_{k/I})\\
           &=: h_{q_{k/i}}(x_i).
 \end{align}
Given the measurement function $h_{q_{k/i}}(\cdot)$, we can now derive the corresponding Jacobian for the measured relative attitude by differentiating $q_{k/i}$ by the state variables $\delta q_{i/I}$ and $\delta q_{k/I}$, as
\begin{align}
  H_{q_{k/i}, \delta q_{i/I}} &= \left.\frac{\partial q_{k/i}}{\partial \delta q_{i/I}}\right|_{x_i = \hat x_i}\\
                              &= \left.\frac{\partial ((\hat q_{i/I}\delta q_{i/I})^*(\hat q_{k/I}\delta q_{k/I}))}{\partial \delta q_{i/I}}\right|_{x_i = \hat x_i}\\
                              &= \rqm{\hat q_{k/I}}\rqm{\hat q_{i/I}^*}\frac{\partial (\delta q_{i/I}^*)}{\partial \delta q_{i/I}}\\
                              &= \rqm{\hat q_{k/I}}\rqm{\hat q_{i/I}^*}I^* \in \mathbb{R}^{4\times 4},
\end{align}
and,
\begin{align}
        H_{q_{k/i}, \delta q_{k/I}} &= \left.\frac{\partial q_{k/i}}{\partial \delta q_{k/I}} \right|_{x_i = \hat x_i}\\
        &= \frac{\partial ((\hat q_{i/I}\delta q_{i/I})^*(\hat q_{k/I}\delta q_{k/I}))}{\partial \delta q_{k/I}}\\
        &= \lqm{\hat q_{i/I}^*}\lqm{\hat q_{k/I}}\frac{\partial (\delta q_{k/I})}{\partial \delta q_{k/I}}\\
        &= \lqm{\hat q_{i/I}^*}\lqm{\hat q_{k/I}} \in \mathbb{R}^{4\times 4}.
\end{align}

\subsubsection{Jacobian of the Relative Position Measurement}
As before we will drop the subscript $B$ denoting the body frame and introduce $q_{r, B_k/B_i} := q_{k/i}$ and $q_{d, B_k/B_I} := p_{k/i}$. 
Satellite $i$ measures its position in its body frame relative to the frame of its neighbour $k$.
Denote this measurement $r_{k/i}^i$.
To compute the measurement function $h_{i,r_{ik}}(\cdot)$, we need two more equations, the position quaternion error and the transformation of the position onto the quaternion position. 
The former is computed as,
\begin{align}
    \dq{\delta q}_{B/I} &= \dq{\hat q}_{B/I}^* \dq{q}_{B/I} \\
  \dq{q}_{B/I} &= \dq{\hat q}_{B/I} \dq{\delta q}_{B/I} \\
                        &= (\hat q_{B/I} \delta q_{B/I}) + \epsilon(\hat q_{B/I} \delta p_{B/I} + \hat p_{B/I} \delta q_{B/I}) \\
                        &= q_{B/I} + \epsilon p_{B/I}
    \end{align}
Hence, $p_{B/I} = \hat q_{B/I} \delta p_{B/I} + \hat p_{B/I} \delta q_{B/I}$.
The latter quantity is simply $r_{B/I}^B = \begin{bmatrix} 0 & \overline{r_{B/I}^B}^\intercal\end{bmatrix}^\intercal = 2p_{B/I}q_{B/I}$. 
We can thus compute the measurement function as,
\begin{align}
        r_{k/I}^I &= r_{i/I}^I + q_{i/I}r_{k/i}^iq_{i/I}^*\\
        r_{k/i}^i &=  q_{i/I}^*(r_{k/I}^I - r_{i/I}^I)q_{i/I}\\
                 &=  q_{i/I}^*(2p_{k/I}q_{k/I}^* - 2p_{i/I}q_{i/I}^*)q_{i/I}\\
                 &=  2(\delta q_{i/I}^*\hat q_{i/I}^*)((\hat p_{k/I}\delta q_{k/I})(\delta q_{k/I}^*\hat q_{k/I}^*))(\hat q_{i/I} \delta q_{i/I})\\
                 & \quad
                 + 2(\delta q_{i/I}^*\hat q_{i/I}^*)((\hat q_{k/I}\delta p_{k/I})(\delta q_{k/I}^*\hat q_{k/I}^*))(\hat q_{i/I} \delta q_{i/I})\\
                 &\quad
                 - 2(\delta q_{i/I}^*\hat q_{i/I}^*)((\hat p_{i/I}\delta q_{i/I})(\delta q_{i/I}^*\hat q_{i/I}^*))(\hat q_{i/I} \delta q_{i/I})\\
                 &\quad
                 -2(\delta q_{i/I}^*\hat q_{i/I}^*)((\hat q_{i/I}\delta p_{i/I})(\delta q_{i/I}^*\hat q_{i/I}^*))(\hat q_{i/I} \delta q_{i/I})\\
                 &= h_{_ir_{ik}}(x_i).
\end{align}
Next, we compute the corresponding Jacobian by differentiating $r_{k/i}^i$ with respect to $\delta q_{i/I}$, $\delta q_{k/I}$, $\delta p_{i/I}$ and $\delta p_{k/I}$. 
First, we can compute
\begin{align}
 & H_{r_{k/i}^i, \delta p_{k/I}} = \left.\frac{\partial r_{k/i}^i}{\partial \delta p_{k/I}}\right|_{x_i = \hat x_i}  \\
                                &= \quad \left.2\frac{\partial (\delta q_{i/I}^*\hat q_{i/I}^*)((\hat q_{k/I}\delta p_{k/I})(\delta q_{k/I}^*\hat q_{k/I}^*))(\hat q_{i/I} \delta q_{i/I})}{\partial \delta p_{k/I}}\right|_{x_i = \hat x_i} \\
                                &= \quad 2\rqm{\hat q_{i/I}^*}\lqm{\hat q_{i/I}}\rqm{\hat q_{k/I}^*}\lqm{\hat q_{k/I}} \in \mathbb{R}^{4\times 4} ,
\end{align}
followed by,
\begin{align}
            &H_{r_{k/i}^i, \delta q_{k/I}} = \left.\frac{\partial r_{k/i}^i}{\partial \delta q_{k/I}}\right|_{x_i = \hat x_i} \\
            &=  2\frac{\partial ((\delta q_{i/I}^*\hat q_{i/I}^*)((\hat p_{k/I}\delta q_{k/I})(\delta q_{k/I}^*\hat q_{k/I}^*))(\hat q_{i/I} \delta q_{i/I})}{\partial \delta q_{k/I}}  + \left.2\frac{\partial (\delta q_{i/I}^*\hat q_{i/I}^*)((\hat q_{k/I}\delta p_{k/I})(\delta q_{k/I}^*\hat q_{k/I}^*))(\hat q_{i/I} \delta q_{i/I}))}{\partial \delta q_{k/I}}\right|_{x_i = \hat x_i}\\
            &=  2\rqm{\hat q_{i/I}}\lqm{\hat q_{i/I}^*}\left(\rqm{\hat q_{k/I}^*}\lqm{\hat p_{k/I}} + \lqm{\hat p_{k/I}}\rqm{\hat q_{k/I}^*}I^*\right) \in \mathbb{R}^{4 \times 4}.\label{eq:sensitivity}
\end{align}
Next we have,
\begin{align}
&            H_{r_{k/i}^i, \delta q_{i/I}} = \left.\frac{\partial r_{k/i}^i}{\partial \delta q_{i/I}}\right|_{x_i = \hat x_i} \\
            &=  2\frac{\partial ((\delta q_{i/I}^*\hat q_{i/I}^*)((\hat p_{k/I}\delta q_{k/I})(\delta q_{k/I}^*\hat q_{k/I}^*))(\hat q_{i/I} \delta q_{i/I}))}{\partial \delta q_{i/I}} \left.
                 - 2\frac{\partial ((\delta q_{i/I}^*\hat q_{i/I}^*)((\hat p_{i/I}\delta q_{i/I})(\delta q_{i/I}^*\hat q_{i/I}^*))(\hat q_{i/I} \delta q_{i/I}))}{\partial \delta q_{i/I}}\right|_{x_i = \hat x_i}\\
                 &=  2\rqm{\hat q_{i/I}}\rqm{\hat p_{k/I}\hat q_{k/I}^*}\rqm{\hat q_{i/I}^*}I^* + 2\lqm{\hat q_{i/I}^*}\lqm{\hat p_{k/I} \hat q_{k/I}^*}\lqm{\hat q_{i/I}}\\
                 & \quad 
                 - 2 \rqm{\hat q_{i/I}}\rqm{\hat p_{i/I}\hat q_{i/I}^*}\rqm{\hat q_{i/I}^*}I^*
                 -2 \lqm{\hat q_{i/I}^*} \lqm{\hat p_{i/I} \hat q_{i/I}^*} \lqm{\hat q_{i/I}}\\
                 & \quad
                 - 2 \lqm{\hat q_{i/I}^*}\rqm{\hat q_{i/I}}\rqm{\hat q_{i/I}^*}\lqm{\hat p_{i/I}} - 2 \lqm{\hat q_{i/I}^*}\rqm{\hat q_{i/I}}\lqm{\hat p_{i/I}}\rqm{q_{i/I}^*}I^*\in \mathbb{R}^{4 \times 4},
\end{align}
and finally,
\begin{align}
            H_{r_{k/i}^i, \delta p_{i/I}} &= \left.\frac{\partial r_{k/i}^i}{\partial \delta p_{i/I}}\right|_{x_i = \hat x_i} \\
            &=  \left.-2\frac{\partial ((\delta q_{i/I}^*\hat q_{i/I}^*)(\hat q_{i/I}\delta p_{i/I}))}{\partial \delta p_{i/I}}\right|_{x_i = \hat x_i}\\
            &=  -2I\in \mathbb{R}^{4 \times 4}.
\end{align}
This concludes the computation of the Jacobians of the individual measurement functions in Eq.~\eqref{eq:measurement_vec_distr}.

\subsection{Full Jacobian for Hard Consensus}
\label{sec:full-jacobian-hard}

In the hard consensus step described in \S\ref{sec:descr-hard-cons}, the measurement Jacobian that satellite $i$ prepares for its neighbour $k$ is given by,
\begin{align}
        \overline{H_{i, k}} &= \begin{bmatrix}
            \overline{H^{n_{k, 1}}_{i, n_{k, 1}}} &  \mathbf{0}_{6 \times 6} & \dots & \overline{H^i_{i, n_{k, 1}}} & \mathbf{0}_{6 \times 6} & \dots & \mathbf{0}_{6 \times 6} & \mathbf{0}_{6 \times 6}\\
            \vdots & \vdots & \ddots & \vdots & \vdots & \ddots & \vdots & \vdots\\
            \mathbf{0}_{6 \times 6} & \mathbf{0}_{6\times 6} & \dots & \mathbf{I}_{6\times 6} & \mathbf{0}_{6 \times 6} & \dots & \mathbf{0}_{6 \times 6} & \mathbf{0}_{6 \times 6}\\
            \vdots & \vdots & \ddots & \vdots & \vdots & \ddots & \vdots & \vdots\\
             \mathbf{0}_{6 \times 6} & \mathbf{0}_{6\times 6} & \dots & \overline{H^i_{i, k_{i, k_k}}} & \mathbf{0}_{6 \times 6} & \dots & \overline{H^{k_k}_{i, k_{i, k_k}}} &  \mathbf{0}_{6 \times 6} 
            \end{bmatrix}\in \mathbb{R}^{6(k_k+1)\times 12(k_k+1)}\label{eq:H_hard_1}\\
        H_{i, k} &= \begin{bmatrix}
            H^{n_{k, 1}}_{i, n_{k, 1}} &  \mathbf{0}_{8 \times 6} & \dots & H^i_{i, n_{k, 1}} & \mathbf{0}_{8 \times 6} & \dots & \mathbf{0}_{8 \times 8} & \mathbf{0}_{8 \times 6}\\
            \vdots & \vdots & \ddots & \vdots & \vdots & \ddots & \vdots & \vdots\\
            \mathbf{0}_{8 \times 8} & \mathbf{0}_{8\times 6} & \dots & \mathbf{I}_{8\times 8} & \mathbf{0}_{8 \times 6} & \dots & \mathbf{0}_{8 \times 6} & \mathbf{0}_{8 \times 6}\\
            \vdots & \vdots & \ddots & \vdots & \vdots & \ddots & \vdots & \vdots\\
             \mathbf{0}_{8 \times 8} & \mathbf{0}_{8\times 6} & \dots & H^i_{i, k_{i, k_k}} & \mathbf{0}_{8 \times 6} & \dots & H^{k_k}_{i, k_{i, k_k}} &  \mathbf{0}_{8 \times 6} 
            \end{bmatrix}\in \mathbb{R}^{8(k_k+1)\times 14(k_k+1)}\label{eq:H_hard_2}.
\end{align}

\if\arxivyes1
   \section{Appendix}
\label{sec:appendix-1}

\subsection{DQ-MEKF}
\label{sec:dq-mekf}

In this appendix, we outline the information form of the single satellite DQ-MEKF as presented in~\cite{7171823}, which will be necessary to derive the distributed version in \S\ref{sec:distr-dual-quat}.

\subsubsection{Time Update}
By following \cite{7171823}, the state and process noise of the extended dual quaternion Kalman filter are
\begin{equation}\label{eq:single_sat_state}
    x_{16} = \begin{bmatrix}\dq{\delta q}_{B/I} \\ \dq{b_{\omega}} \end{bmatrix} \in \mathbb{R}^{16} \quad \text{and}\quad w_{16} = \begin{bmatrix}\dq{\eta}_{\omega} \\ \dq{\eta}_{b_{\omega}} \end{bmatrix} \in \mathbb{R}^{16},
\end{equation}
where $\dq{\delta q}_{B/I}$ is defined as the \emph{error dual quaternion}, and $\dq{b_{\omega}}$ is the \emph{dual bias} of the velocity. 
 The bias acting on the angular velocity measurement is denoted by $\overline{b_\omega} \in \mathbb{R}^3$, and the bias acting on the linear velocity measurement is denoted by $\overline{b_v} \in \mathbb{R}^3$. 
Explicitly, these quantities are,
\begin{equation}
    \dq{\delta q}_{B/I} = \dq{\hat q}_{B/I}^*\dq{q}_{B/I} \in \mathbb{H}^u_d \quad \text{and}\quad \dq{b}_{\omega} = b_{\omega} + \epsilon b_v \in \mathbb{H}^v_d.
\end{equation}

Eq.~\eqref{eq:derivation_kinematics} is used to propagate the pose of the satellite. 
It is important to point out the state as defined in \eqref{eq:single_sat_state} is not propagated, as that state is the \emph{error} dual quaternion. 
Rather, the pose $\dq{q}_{B/I}$ is propagated, since ultimately, this is the quantity which we wish to estimate.
The update equation is therefore,
\begin{subequations}\label{eq:single_sat_prop}
  \begin{align}
    \dot{\dq{q}}_{B/I} &= \frac{1}{2}\dq{\omega}_{B/I}^B\dq{q}_{B/I}\\
    \dot{\hat{\bf{q}}}_{B/I} &= \frac{1}{2}\dq{\hat \omega}_{B/I}^B\dq{\hat q}_{B/I}.
  \end{align}
\end{subequations}
We use the velocity measurement model,
  \begin{align}
    \dq{\omega}_{B/I, m}^B &= \dq{\omega}_{B/I}^B+\dq{b}_{\omega} + \dq{\eta}_{\omega}\label{eq:velocity_measurement_model_1}\\ 
    \dq{\hat \omega}_{B/I}^B &= \dq{\omega}_{B/I, m}^B - \dq{\hat b}_{\omega},\label{eq:velocity_measurement_model_2}
  \end{align}
where $\dq{\omega}_{B/I, m}^B$ is the measured dual velocity in the body frame, $\dq{\omega}_{B/I}^B$ is the true dual velocity in the body frame, and $\dq{b}_{\omega}$ is the dual bias. 
We denote $\dq{\eta}_{\omega}$ as the Gaussian noise acting on the velocity measurement, with $\text{E}[\overline{\dq{\eta}}_{\omega}] = \mathbf{0}_{6\times 1}$ and $\Cov[\overline{\dq{\eta}}_{\omega}] = \overline Q_{\dq{\omega}}$. 
As is common with models of rate integration gyros, the bias is assumed to be driven by another zero mean Gaussian noise process $\dq{\eta}_{b_{\omega}}$ with $\text{E}[\overline{\dq{\eta}}_{b_{\omega}}] = \mathbf{0}_{6\times 1}$ and $\Cov[\overline{\dq{\eta}}_{b_{\omega}}] = \overline Q_{\dq{b_{\omega}}}$.
The covariance matrices $\overline Q_{\dq{\omega}}$ and $\overline Q_{\dq{b_{\omega}}}$ are given by, 
\begin{equation}\label{eq:single_sat_Q}
  \overline Q_{\dq{\omega}} = 
  \begin{bmatrix}
    \overline Q_{\omega} & \mathbf{0}_{3\times 3} \\ 
    \mathbf{0}_{3\times 3} &\overline Q_v 
  \end{bmatrix},\qquad 
  \overline Q_{\dq{b_{\omega}}} = 
  \begin{bmatrix}
    \overline Q_{b_{\omega}} & \mathbf{0}_{3\times 3} \\ 
    \mathbf{0}_{3\times 3} &\overline Q_{b_v},
  \end{bmatrix}
\end{equation}
and finally, the dual bias satisfies
\begin{equation}\label{eq:bias_propagation_equation}
    \dot{\dq{b}}_{\omega} = \dq{\eta}_{b_{\omega}}.
\end{equation}
Following \cite{7171823}, the error quaternion dynamics can be defined to propagate the (reduced) state variables\footnote{As discussed in~\cite{optimal_estimation_of_dynamic_syxstems}, the unit quaternion constraint allows a dimensional reduction of the state variable, as the scalar parts of the dual quaternion are uniquely recovered by their corresponding vector parts},
\begin{align}
    x = \begin{bmatrix}\overline{\dq{\delta q}_{B/I}} \\ \overline{\dq{b_{\omega}}} \end{bmatrix} \in \mathbb{R}^{12} \quad \text{and}\quad w = \begin{bmatrix}\overline{\dq{\eta_{\omega}}} \\ \overline{\dq{\eta_{b_{\omega}}}} \end{bmatrix} \in \mathbb{R}^{12}.
\end{align}
The corresponding functions $f$ and $g$ for the dynamics~\eqref{eq:def_kalman_state_eq} are,
\begin{align}
      f(x(t)) &= 
                        \begin{bmatrix}
                          \overline{\frac{1}{2}\dq{\delta q}_{B/I}\dq{\hat \omega}_{B/I}^B - \frac{1}{2}\dq{\hat \omega}_{B/I}^B\dq{\delta q}_{B/I} + \frac{1}{2}\dq{\delta q}_{B/I}\dq{\hat b_{\omega}} - \frac{1}{2}\dq{\delta q}_{B/I}\dq{b_{\omega}} - \frac{1}{2}\dq{\delta q}_{B/I}\dq{\eta_{\omega}}}\\ 
                          \mathbf{0}_{6\times 1}
                        \end{bmatrix}\label{eq:f12_1}\\
    g(x(t)) &= 
                                 \begin{bmatrix}
                                   -\frac{1}{2}[\dq{\tilde{\delta q}}_{B/I}] & \mathbf{0}_{6 \times 6}
                                   \\ \mathbf{0}_{6\times 6} & \mathbf{I}_{6 \times 6}
                                 \end{bmatrix}\label{eq:g12_1}
\end{align}
and their corresponding Jacobians are,
\begin{align}
    F(t) &= \frac{\partial f(x)}{\partial x} = 
                         \begin{bmatrix}
                           -\overline{\dq{\hat \omega}_{B/I}^B}^{\times} & -\frac{1}{2}\mathbf{I}_{6\times 6} \\ 
                           \mathbf{0}_{6\times 6} & \mathbf{0}_{6\times 6}
                         \end{bmatrix},~
    G(t) = g(x)|_{x = \hat x} = 
                         \begin{bmatrix}
                           -\frac{1}{2}I_{6 \times 6} & \mathbf{0}_{6 \times 6} \\ 
                           \mathbf{0}_{6 \times 6} & \mathbf{I}_{6\times 6}
                         \end{bmatrix}.\label{Single_Sat_F_G_Mat}
\end{align}
With the reduced states, the process noise matrix $Q$ becomes
\begin{equation}
  Q = \mathrm{blkdiag} \left\{ \overline{Q}_{\dq{\omega}},\overline{Q}_{\dq{{b_\omega}}}  \right\}.
\end{equation}


\subsubsection{Measurement Update}
For the measurement update the filter needs a measurement of the pose, i.e.~a measurement of the attitude and position in reference to an inertial frame. 
We denote that quantity as $\dq{q}_{B/I, m}$, leading to the measurement model 
\begin{align}\label{eq:single_satellite_h}
    z_{m}(t_k) &= \overline{\dq{\hat{q}}^{-}_{B/I}(t_k)^*\dq{q}_{B/I, m}}(t_k)\\
    & = h(x(t_k)) + v(t_k) = \overline{\dq{\delta q}_{B/I}(t_k)} + v(t_k),
\end{align}
where $v_6$ is the Gaussian noise on the measurement with $E[v] = \mathbf{0}$ and $\Cov[v] = R$. 
The Jacobian $H$ of the measurement function $h$ (also called the \emph{sensitivity}) is then,
\begin{equation}\label{eq:single_satellite_H}
    h(x(t_k)) = \overline{\dq{\delta q}_{B/I}(t_k)}\quad \rightarrow \quad H = 
    \begin{bmatrix}
      \mathbf{I}_{6\times 6} & \mathbf{0}_{6\times 6}
    \end{bmatrix}.
\end{equation}
Next, using the quantities $u$, $U$ and $M$ in Eqs.~\eqref{eq:y_S_inf_form_def} and \eqref{eq:M_inf_form_def}, the optimal Kalman update is
  \begin{align}\label{eq:single_satellite_kalman_update}
    \Delta^* x(t_k) &= 
    \begin{bmatrix}
      \Delta^* \overline{\dq{\delta \hat{q}}_{B/I}(t_k)} \\ 
      \Delta^* \overline{\dq{\delta \hat{b}_\omega}(t_k)} 
    \end{bmatrix} \\
    &= M(t_k)[u(t_k) - U(t_k) \hat x(t_k)] \\
    &= M(t_k)H^\intercal R^{-1}\text{ } \overline{\dq{\hat{q}}^{-}_{B/I}(t_k)^*\dq{q}_{B/I, m}(t_k)}.
    \end{align}
This yields the unit dual quaternion,
\begin{equation}\label{eq:extend_to_R8}
\begin{aligned}
&\Delta^* \delta \hat q_{B/I, r} = \left(\sqrt{1-||\Delta^*\overline{\delta \hat{q}_{B/I, r}(t_k)}||^2}, \Delta^*\overline{\delta \hat{q}_{B/I, r}}(t_k)\right) \\
&\Delta^* \delta \hat q_{B/I, d}= \left(\frac{-\Delta^*\overline{\delta \hat{q}_{B/I, r}(t_k)}^T\Delta^*\overline{\delta \hat{q}_{B/I, d}(t_k)}}{\sqrt{1-||\Delta^*\overline{\delta \hat{q}_{B/I, r}(t_k)}||^2}}, \Delta^*\overline{\delta \hat{q}_{B/I, d}}(t_k)\right),
\end{aligned}
\end{equation}
and the measurement update is concluded by computing the estimate
\begin{subequations}\label{eq:single_Sat_meas_update_estimated}
\begin{align}
    \dq{\hat{q}}^+_{B/I}(t_k) &= \dq{\hat{q}}^-_{B/I}(t_k)\Delta^* \dq{\delta \hat q}_{B/I}(t_k)\\
    \overline{\dq{\hat{b}^+_\omega}(t_k)} &= \overline{\dq{\hat{b}^-_\omega}(t_k)} + \Delta^* \overline{\dq{\delta \hat{b}}_\omega(t_k)}.
\end{align}
\end{subequations}

In the following subsections, we summarize two scenarios where different measurements are available.

\subsubsection{Scenario 1: Pose measurements only}\label{sec:single_KF_pose_meas_only}
In the above derivation, it was assumed that linear and angular velocity measurements are accessible. 
In some applications, only pose measurements might be accessible. 
For example, in the distributed case discussed in \S\ref{sec:distr-dual-quat}, a satellite outfitted with inexpensive cameras running a pose detection algorithm~\cite{Kreiss_2019_CVPR,kreiss2021openpifpaf}.
This scenario is easily covered by changing the velocity measurement model.
If the velocity measurement is not available $\dq{\omega}^B_{B/I, m}$ and $\dq{\eta}_\omega$ are set to zero, effectively leaving us with
\begin{subequations}\label{eq:pose_meas_only_velocity_model}
  \begin{align}
    \dq{\omega}^B_{B/I} &= -\dq{b}_\omega \\
    \dq{\hat\omega}^B_{B/I} &= -\dq{\hat b}_\omega.
  \end{align}
\end{subequations}
Simply setting $\overline{Q_{\dq{\omega}}}$ to zero covers this case. 
The estimated velocity is then just the negative of the bias as shown in Eq.~\eqref{eq:pose_meas_only_velocity_model}. 

\subsubsection{Pose, angular velocity, and linear acceleration measurements}
Inertial measurements units (IMUs) are common in satellite hardware, and so the case where pose, angular velocity, and linear acceleration measurements are available is covered in this section.

Those sensors provide the angular velocity and the linear acceleration. 
Unfortunately, when using this setting, the duality of attitude and position in the dual quaternion has to be broken.
Following \cite{7171823}, the linear acceleration model is defined as
\begin{equation}\label{eq:lin_vel_model}
    n_{A/I, m}^B = n_{A/I}^B + b_n + \eta_n \in \mathbb{H}^v,
\end{equation}
where $n_{A/I}^B$ is the force at the location of the IMU with respect to the inertial frame expressed in body frame. 
The bias of the IMU is denoted by $b_n\in\mathbb{H}^v$, and $\eta_n \in \mathbb{H}^v$ is the noise acting on the true acceleration with $E[\overline{\eta_n}] = \mathbf{0}_{3\times 1}$ and $\Cov[\overline{\eta_n}] = \overline{Q_n} \in \mathbb{R}^{3\times 3}$. 
The bias of the IMU is assumed to be driven by zero mean Gaussian process
\begin{equation}
    \dot b_n = \eta_{b_n} \in \mathbb{H}^v,
\end{equation}
with $E[\overline{\eta_{b_n}}] = \mathbf{0}_{3\times 1}$ and $\Cov[\overline{\eta_{b_n}}] = \overline{Q_{b_n}} \in \mathbb{R}^{3\times 3}$. 
The estimate of the acceleration is given through Eq.~\eqref{eq:lin_vel_model}, where the noise in the given timestep is assumed to be zero:
\begin{equation}\label{eq:lin_acc_estimation}
    \hat n_{A/I}^B = n_{A/I, m}^B - \hat b_n.
\end{equation}
The state of the DQ-MEKF needs to be extended to incorporate the IMU bias:
\begin{align}
  \label{eq:single_sat_state_imu}
  \begin{split}
    x &= \begin{bmatrix}
      \dq{\delta q}_{B/I}^T & 
      \dq{b_{\omega}}^T& 
      b_n^T 
    \end{bmatrix}^T \in \mathbb{R}^{20}\\
    w &= \begin{bmatrix}
      \dq{\eta}_{\omega}^T & 
      \dq{\eta}_{b_{\omega}}^T & 
      \eta_{b_n}^T 
    \end{bmatrix}^T \in \mathbb{R}^{20}.
  \end{split}
\end{align}
As there are no linear velocity measurements, the estimate of the linear velocity becomes
\begin{equation}\label{eq:lin_acc_kf_lin_vel}
    \hat v_{B/I}^B = -\hat b_v \qquad v_{B/I}^B = - b_v.
\end{equation}
Hence, $Q_v = \mathbf{0}_{3\times 3}$.
By following \cite{7171823} and neglecting gravity, the time derivative of \eqref{eq:lin_acc_kf_lin_vel} and thus the filter dynamics use the Jacobians,
\begin{align}
    F &= \begin{bmatrix}-\overline{\hat \omega_{B/I}^B}^\times & \mathbf{0}_{3\times 3} & -\frac{1}{2}\mathbf{I}_{3\times 3} & \mathbf{0}_{3\times 3} & \mathbf{0}_{3\times 3} \\
    -\overline{\hat v_{B/I}^B}^\times & -\overline{\hat \omega_{B/I}^B}^\times & \mathbf{0}_{3\times 3} & -\frac{1}{2}\mathbf{I}_{3\times 3} & \mathbf{0}_{3\times 3} \\
    \mathbf{0}_{3\times 3} & \mathbf{0}_{3\times 3} & \mathbf{0}_{3\times 3} & \mathbf{0}_{3\times 3} & \mathbf{0}_{3\times 3}\\
    \mathbf{0}_{3\times 3} & \mathbf{0}_{3\times 3} & -\overline{\hat b_v}^\times + \overline{\hat \omega_{B/I}^B \times r_{A/B}^B}^\times + \overline{\hat \omega_{B/I}^B}^\times \overline{r_{A/B}^B}^\times & -\overline{\hat \omega_{B/I}^B}^\times & \mathbf{I}_{3\times 3} \\
    \mathbf{0}_{3\times 3} & \mathbf{0}_{3\times 3} & \mathbf{0}_{3\times 3} & \mathbf{0}_{3\times 3} & \mathbf{0}_{3\times 3}
    \end{bmatrix}\label{eq:F_linear_acc}\\
    G
                    & = \mathrm{blkdiag} \left\{ -\frac{1}{2} \mathbf{I}_{6\times6}, \mathbf{I}_{9\times 9} \right\}
                    + 
                      \begin{bmatrix}
                        \mathbf{0}_{9\times 3} & | \\ 
                        -\overline{\hat b_v}^\times + \overline{\hat \omega_{B/I}^B \times r_{A/B}^B}^\times + \overline{\hat \omega_{B/I}^B}^\times \overline{r_{A/B}^B}^\times & \mathbf{0}_{15\times 12}\\
                        \mathbf{0}_{3\times 3} & |
                      \end{bmatrix}\label{eq:G_linear_acc},
\end{align}

where the states are reduced by neglecting the scalar parts of the dual quaternion.
The translation vector from the origin of the body frame to the location of the IMU expressed in the body frame is denoted by $r_{A/B}^B \in \mathbb{H}^v$. 
As before, the states where reduced.
\begin{align}
  x &= 
  \begin{bmatrix}
    \overline{\dq{\delta q}_{B/I}}^T & 
    \overline{\dq{b_{\omega}}}^T & 
    \overline{b_n}^T 
  \end{bmatrix}^T \in \mathbb{R}^{15} \\
  w &= \begin{bmatrix}
    \overline{\dq{\eta_{\omega}}}^T & 
    \overline{\dq{\eta_{b_{\omega}}}}^T & 
    \overline{\eta_{b_n}}^T 
  \end{bmatrix}^T \in \mathbb{R}^{15}
\end{align}
Following~\cite{7171823}, the linear velocity bias is propagated using the approximation,
\begin{align}
    \dot{\hat b}_v \approx -\hat \omega_{B/I}^B \times \hat b_v - \hat n_{A/I}^B + \hat \omega_{B/I}^B\times\left(\hat \omega_{B/I}^B\times r_{A/B}^B\right).\label{eq:lin_bias_prop}
\end{align}
For the time update of this filter variant, the acceleration using Eq.~\eqref{eq:lin_acc_estimation} is computed first, followed by the linear velocity bias using Eq.~\eqref{eq:lin_bias_prop}.
The linear velocity is computed using Eq.~\eqref{eq:lin_acc_kf_lin_vel}, and the angular velocity is computed using the real part of Eq.~\eqref{eq:velocity_measurement_model_2}. 
This concludes the calculation of the estimated dual velocity. 
Then, using Eq.~\eqref{eq:single_sat_prop} $\dq{\hat q_{B/I}}$ is computed.

The measurement update is computed the measurement matrix,
\begin{equation}
    H = \begin{bmatrix} \mathbf{I}_{6\times 6} & \mathbf{0}_{6\times 6} & \mathbf{0}_{6\times 3} \end{bmatrix}
\end{equation}
and the optimal Kalman update is therefore,
\begin{align}\label{eq:single_satellite_kalman_update_lin_acc}
    \Delta^* x_{15}(t_k) &= 
    \begin{bmatrix}
      \Delta^* \overline{\dq{\delta \hat{q}}_{B/I}(t_k)} \\ 
      \Delta^* \overline{\dq{\delta \hat{b}_\omega}(t_k)} \\ 
      \Delta^* \overline{\delta \hat{b}_n(t_k)}
    \end{bmatrix} 
    = M(t_k)H^\intercal R^{-1}\text{ } \overline{\dq{\hat{q}}^{-}_{B/I}(t_k)^*\dq{q}_{B/I, m}(t_k)}.
\end{align}
The estimates are then updated using Eqs.~\eqref{eq:extend_to_R8} and \eqref{eq:single_Sat_meas_update_estimated}. 
Finally, the linear acceleration bias is then updated by 
\begin{equation}
\overline{\hat b_n^+(t_k)} = \overline{\hat b_n^-(t_k)} + \Delta^*\overline{\delta \hat b_n(t_k)}.
\end{equation}

\fi

\section*{Acknowledgments}
This research was funded by the Swiss National Science Foundation under the NCCR Automation (180545) and the Research Council of Norway under the Centre for Space Sensors \& Systems (309835).


\begin{thebibliography}{40}
\newcommand{\enquote}[1]{``#1''}
\providecommand{\natexlab}[1]{#1}
\providecommand{\url}[1]{\texttt{#1}}
\providecommand{\urlprefix}{URL }
\expandafter\ifx\csname urlstyle\endcsname\relax
  \providecommand{\doi}[1]{\discretionary{}{}{}https://doi.org/#1}\else
  \providecommand{\doi}[1]{\discretionary{}{}{}\urlstyle{rm}\url{https://doi.org/#1}}\fi

\bibitem[{Johnson(2010)}]{orbital_debris}
Johnson, N., \enquote{Orbital debris: The growing threat to space operations,}
  \emph{Advances in the Astronautical Sciences}, Vol. 137, 2010, pp. 3--11.

\bibitem[{Bandyopadhyay and Quadrelli(2017)}]{bandyopadhyay2017optimal}
Bandyopadhyay, S., and Quadrelli, M., \enquote{Optimal Transport Based Control
  of Granular Imaging System in Space,} \emph{9th International Workshop on
  Satellite Constellations and Formation Flying}, 2017, pp. 1--11.
\newblock \doi{10.2514/1.G006001}.

\bibitem[{Danzmann and {the LISA Study Team}(1996)}]{danzmann1996lisa}
Danzmann, K., and {the LISA Study Team}, \enquote{LISA: laser interferometer
  space antenna for gravitational wave measurements,} \emph{Classical and
  Quantum Gravity}, Vol.~13, No. 11A, 1996, p. A247.
\newblock \doi{10.2514/6.1995-829}.

\bibitem[{Amaro-Seoane et~al.(2023)Amaro-Seoane, Andrews, Arca~Sedda, Askar,
  Baghi, Balasov, Bartos, Bavera, Bellovary, Berry et~al.}]{amaro2017laser}
Amaro-Seoane, P., Andrews, J., Arca~Sedda, M., Askar, A., Baghi, Q., Balasov,
  R., Bartos, I., Bavera, S.~S., Bellovary, J., Berry, C.~P., et~al.,
  \enquote{Astrophysics with the laser interferometer space antenna,}
  \emph{Living Reviews in Relativity}, Vol.~26, No.~1, 2023, p.~2.
\newblock \doi{10.1007/s41114-022-00041-y}.

\bibitem[{Markley(2003)}]{Attitude_Error_Representations}
Markley, L., \enquote{Attitude Error Representations for Kalman Filtering,}
  \emph{Journal of Guidance Control and Dynamics}, Vol.~26, 2003, pp. 311--317.
\newblock \doi{10.2514/2.5048}.

\bibitem[{Farrell(1970)}]{FARRELL1970419}
Farrell, J., \enquote{Attitude determination by {Kalman} filtering,}
  \emph{Automatica}, Vol.~6, No.~3, 1970, pp. 419--430.
\newblock \doi{10.1016/0005-1098(70)90057-9}.

\bibitem[{Crassidis et~al.(2007)Crassidis, Markley, and
  Cheng}]{Survey_of_Nonlinear_Attitude_Estimation_Methods}
Crassidis, J., Markley, L., and Cheng, Y., \enquote{Survey of Nonlinear
  Attitude Estimation Methods,} \emph{Journal of Guidance Control and
  Dynamics}, Vol.~30, 2007, pp. 12--28.
\newblock \doi{10.2514/1.22452}.

\bibitem[{{Filipe} et~al.(2015){Filipe}, {Kontitsis}, and {Tsiotras}}]{7171823}
{Filipe}, N., {Kontitsis}, M., and {Tsiotras}, P., \enquote{Extended Kalman
  Filter for spacecraft pose estimation using dual quaternions,} \emph{2015
  American Control Conference (ACC)}, 2015, pp. 3187--3192.
\newblock \doi{10.1109/ACC.2015.7171823}.

\bibitem[{{Michieletto} and {Cenedese}(2019)}]{8795869}
{Michieletto}, G., and {Cenedese}, A., \enquote{Formation Control for Fully
  Actuated Systems: a Quaternion-based Bearing Rigidity Approach,} \emph{2019
  18th European Control Conference (ECC)}, 2019, pp. 107--112.
\newblock \doi{10.23919/ECC.2019.8795869}.

\bibitem[{Lee and Mesbahi(2017)}]{lee2017constrained}
Lee, U., and Mesbahi, M., \enquote{Constrained autonomous precision landing via
  dual quaternions and model predictive control,} \emph{Journal of Guidance,
  Control, and Dynamics}, Vol.~40, No.~2, 2017, pp. 292--308.
\newblock \doi{10.2514/1.G001879}.

\bibitem[{Zinage and Bakolas(2022)}]{koopmandq}
Zinage, V., and Bakolas, E., \enquote{Koopman operator based modeling and
  control of rigid body motion represented by dual quaternions,} \emph{2022
  American Control Conference (ACC)}, 2022, pp. 3997--4002.
\newblock \doi{10.23919/ACC53348.2022.9867584}.

\bibitem[{{Zu} et~al.(2014){Zu}, {Sun}, and {Dai}}]{7039403}
{Zu}, Y., {Sun}, C., and {Dai}, R., \enquote{Distributed estimation for spatial
  rigid motion based on dual quaternions,} \emph{53rd IEEE Conference on
  Decision and Control}, 2014, pp. 334--339.
\newblock \doi{10.1109/CDC.2014.7039403}.

\bibitem[{Dong et~al.(2016)Dong, Hu, and Ma}]{DONG201687}
Dong, H., Hu, Q., and Ma, G., \enquote{Dual-quaternion based fault-tolerant
  control for spacecraft formation flying with finite-time convergence,}
  \emph{ISA Transactions}, Vol.~61, 2016, pp. 87--94.
\newblock \doi{10.1016/j.isatra.2015.12.008}.

\bibitem[{Filipe and Tsiotras(2015)}]{filipe2015adaptive}
Filipe, N., and Tsiotras, P., \enquote{Adaptive position and attitude-tracking
  controller for satellite proximity operations using dual quaternions,}
  \emph{Journal of Guidance, Control, and Dynamics}, Vol.~38, No.~4, 2015, pp.
  566--577.
\newblock \doi{10.2514/1.G000054}.

\bibitem[{Yang and Stoll(2019)}]{yang2019adaptive}
Yang, J., and Stoll, E., \enquote{Adaptive sliding mode control for spacecraft
  proximity operations based on dual quaternions,} \emph{Journal of Guidance,
  Control, and Dynamics}, Vol.~42, No.~11, 2019, pp. 2356--2368.
\newblock \doi{10.2514/1.G004435}.

\bibitem[{Reynolds and Mesbahi(2018)}]{coupled_6-DOF_control}
Reynolds, T.~P., and Mesbahi, M., \enquote{Coupled 6-DOF Control for
  Distributed Aerospace Systems,} \emph{2018 IEEE Conference on Decision and
  Control (CDC)}, 2018, pp. 5294--5299.
\newblock \doi{10.1109/CDC.2018.8619618}.

\bibitem[{Crassidis and Junkins(2004)}]{optimal_estimation_of_dynamic_syxstems}
Crassidis, J.~L., and Junkins, J.~L., \emph{{Optimal Estimation of Dynamic
  Systems}}, CRC Press, 2004.
\newblock \doi{10.1201/b11154}.

\bibitem[{{Sabatini}(2006)}]{1643403}
{Sabatini}, A.~M., \enquote{Quaternion-based extended Kalman filter for
  determining orientation by inertial and magnetic sensing,} \emph{IEEE
  Transactions on Biomedical Engineering}, Vol.~53, No.~7, 2006, pp.
  1346--1356.
\newblock \doi{10.1109/TBME.2006.875664}.

\bibitem[{Filipe et~al.(2015)Filipe, Kontitsis, and
  Tsiotras}]{filipe2015extendedjgcd}
Filipe, N., Kontitsis, M., and Tsiotras, P., \enquote{Extended Kalman filter
  for spacecraft pose estimation using dual quaternions,} \emph{Journal of
  Guidance, Control, and Dynamics}, Vol.~38, No.~9, 2015, pp. 1625--1641.
\newblock \doi{10.2514/1.G000977}.

\bibitem[{Alonso et~al.(2001)Alonso, Du, Hughes, Junkins, and
  Crassidis}]{visnav2}
Alonso, R., Du, J.-Y., Hughes, D., Junkins, J.~L., and Crassidis, J.~L.,
  \enquote{Relative navigation for formation flying of spacecraft,} \emph{2001
  Flight Mechanics Symposium}, 2001.

\bibitem[{Alonso et~al.(2000)Alonso, Crassidis, and Junkins}]{alonso2000vision}
Alonso, R., Crassidis, J., and Junkins, J., \enquote{Vision-based relative
  navigation for formation flying of spacecraft,} \emph{AIAA guidance,
  navigation, and control conference and exhibit}, 2000.
\newblock \doi{10.2514/6.2000-4439}.

\bibitem[{Sung et~al.(2022)Sung, Peck, Majji, and Junkins}]{sung2022optical}
Sung, K., Peck, C., Majji, M., and Junkins, J.~L., \enquote{An Optical
  Navigation System for Autonomous Aerospace Systems,} \emph{IEEE Sensors
  Journal}, Vol.~22, No.~17, 2022, pp. 16862--16873.
\newblock \doi{10.1109/JSEN.2022.3190910}.

\bibitem[{Mourikis et~al.(2009)Mourikis, Trawny, Roumeliotis, Johnson, Ansar,
  and Matthies}]{mourikis2009vision}
Mourikis, A.~I., Trawny, N., Roumeliotis, S.~I., Johnson, A.~E., Ansar, A., and
  Matthies, L., \enquote{Vision-aided inertial navigation for spacecraft entry,
  descent, and landing,} \emph{IEEE Transactions on Robotics}, Vol.~25, No.~2,
  2009, pp. 264--280.
\newblock \doi{10.1109/TRO.2009.2012342}.

\bibitem[{Mahmoud and Khalid(2013)}]{dkf_review}
Mahmoud, M., and Khalid, H.~M., \enquote{Distributed Kalman filtering: A
  bibliographic review,} \emph{Control Theory \& Applications, IET}, Vol.~7,
  2013, pp. 483--501.
\newblock \doi{10.1049/iet-cta.2012.0732}.

\bibitem[{Zivan and Choukroun(2018)}]{zivan2018dual}
Zivan, Y., and Choukroun, D., \enquote{Dual quaternion Kalman filters for
  spacecraft relative navigation,} \emph{2018 AIAA Guidance, Navigation, and
  Control Conference}, 2018, p. 1347.
\newblock \doi{10.2514/6.2018-1347}.

\bibitem[{{Olfati-Saber}(2005)}]{1583486}
{Olfati-Saber}, R., \enquote{Distributed Kalman Filter with Embedded Consensus
  Filters,} \emph{Proceedings of the 44th IEEE Conference on Decision and
  Control}, 2005, pp. 8179--8184.
\newblock \doi{10.1109/CDC.2005.1583486}.

\bibitem[{{Olfati-Saber} and {Shamma}(2005)}]{1583238}
{Olfati-Saber}, R., and {Shamma}, J.~S., \enquote{Consensus Filters for Sensor
  Networks and Distributed Sensor Fusion,} \emph{Proceedings of the 44th IEEE
  Conference on Decision and Control}, 2005, pp. 6698--6703.
\newblock \doi{10.1109/CDC.2005.1583238}.

\bibitem[{{Olfati-Saber}(2007)}]{4434303}
{Olfati-Saber}, R., \enquote{Distributed Kalman filtering for sensor networks,}
  \emph{2007 46th IEEE Conference on Decision and Control}, 2007, pp.
  5492--5498.
\newblock \doi{10.1109/CDC.2007.4434303}.

\bibitem[{Ryu and Back(2019)}]{ryu2019distributed}
Ryu, K., and Back, J., \enquote{Distributed Kalman-filtering: Distributed
  optimization viewpoint,} \emph{2019 IEEE 58th Conference on Decision and
  Control (CDC)}, 2019, pp. 2640--2645.
\newblock \doi{10.1109/CDC40024.2019.9029645}.

\bibitem[{Li et~al.(2019)Li, Dong, Li, and Wang}]{LI2019104500}
Li, C., Dong, H., Li, J., and Wang, F., \enquote{Distributed Kalman filtering
  for sensor network with balanced topology,} \emph{Systems \& Control
  Letters}, Vol. 131, 2019, p. 104500.
\newblock \doi{https://doi.org/10.1016/j.sysconle.2019.104500}.

\bibitem[{{Roumeliotis} and {Bekey}(2002)}]{1067998}
{Roumeliotis}, S.~I., and {Bekey}, G.~A., \enquote{Distributed multirobot
  localization,} \emph{IEEE Transactions on Robotics and Automation}, Vol.~18,
  No.~5, 2002, pp. 781--795.
\newblock \doi{10.1109/TRA.2002.803461}.

\bibitem[{Binz et~al.(2023 [online]{\natexlab{a}})Binz, Hudoba~de Badyn,
  Iannelli, and Smith}]{binzcode}
Binz, J., Hudoba~de Badyn, M., Iannelli, A., and Smith, R.~S.,
  \emph{Supplementary software for ``Distributed Dual Quaternion Extended
  Kalman Filtering for Spacecraft Pose Estimation''}, 2023
  [online]{\natexlab{a}}.
\newblock \doi{10.3929/ethz-b-000510769}.

\bibitem[{Binz et~al.(2023 [online]{\natexlab{b}})Binz, Hudoba~de Badyn,
  Iannelli, and Smith}]{arxivversion}
Binz, J., Hudoba~de Badyn, M., Iannelli, A., and Smith, R.~S.,
  \emph{Supplementary reference for ``Distributed Dual Quaternion Extended
  Kalman Filtering for Spacecraft Pose Estimation''}, 2023
  [online]{\natexlab{b}}.
\newblock \doi{10.3929/ethz-b-000619855}.

\bibitem[{{Olfati-Saber}(2009)}]{5399678}
{Olfati-Saber}, R., \enquote{Kalman-Consensus Filter: Optimality, stability,
  and performance,} \emph{Proceedings of the 48h IEEE Conference on Decision
  and Control (CDC) held jointly with 2009 28th Chinese Control Conference},
  2009, pp. 7036--7042.
\newblock \doi{10.1109/CDC.2009.5399678}.

\bibitem[{{Gunnam} et~al.(2002){Gunnam}, {Hughes}, {Junkins}, and
  {Kehtarnavaz}}]{visnav}
{Gunnam}, K.~K., {Hughes}, D.~C., {Junkins}, J.~L., and {Kehtarnavaz}, N.,
  \enquote{A vision-based DSP embedded navigation sensor,} \emph{IEEE Sensors
  Journal}, Vol.~2, No.~5, 2002, pp. 428--442.
\newblock \doi{10.1109/JSEN.2002.806212}.

\bibitem[{Kreiss et~al.(2021)Kreiss, Bertoni, and Alahi}]{kreiss2021openpifpaf}
Kreiss, S., Bertoni, L., and Alahi, A., \enquote{{OpenPifPaf:} Composite fields
  for semantic keypoint detection and spatio-temporal association,} \emph{IEEE
  Transactions on Intelligent Transportation Systems}, Vol.~23, No.~8, 2021,
  pp. 13498--13511.
\newblock \doi{10.1109/TITS.2021.3124981}.

\bibitem[{Kreiss et~al.(2019)Kreiss, Bertoni, and Alahi}]{Kreiss_2019_CVPR}
Kreiss, S., Bertoni, L., and Alahi, A., \enquote{PifPaf: Composite Fields for
  Human Pose Estimation,} \emph{Proceedings of the IEEE/CVF Conference on
  Computer Vision and Pattern Recognition (CVPR)}, 2019, pp. 11977--11986.
\newblock \doi{10.1109/CVPR.2019.01225}.

\bibitem[{Markley et~al.(2007)Markley, Cheng, Crassidis, and
  Oshman}]{quaternion_averaging}
Markley, L., Cheng, Y., Crassidis, J., and Oshman, Y., \enquote{Averaging
  Quaternions,} \emph{Journal of Guidance, Control, and Dynamics}, Vol.~30,
  2007, pp. 1193--1196.
\newblock \doi{10.2514/1.28949}.

\bibitem[{Gonzalez(2010)}]{fibonacci_lattice}
Gonzalez, A., \enquote{Measurement of Areas on a Sphere Using Fibonacci and
  Latitude Longitude Lattices,} \emph{Mathematical Geosciences}, Vol.~42, 2010,
  pp. 49--64.
\newblock \doi{10.1007/s11004-009-9257-x}.

\bibitem[{Binz et~al.(2023 [online]{\natexlab{c}})Binz, Hudoba~de Badyn,
  Iannelli, and Smith}]{ourvideo}
Binz, J., Hudoba~de Badyn, M., Iannelli, A., and Smith, R.~S.,
  \emph{Supplementary video for ``Distributed Dual Quaternion Extended Kalman
  Filtering for Spacecraft Pose Estimation''}, 2023 [online]{\natexlab{c}}.
\newblock \doi{10.3929/ethz-b-000606978}.

\end{thebibliography}
\end{document}